%% file: Article-NonHolonomic_corrected_ARXIV.tex
\documentclass[final]{siamart190516}

\usepackage[colorinlistoftodos]{todonotes}
\usepackage[pagewise]{lineno}

\newcommand\blfootnote[1]{%
  \begingroup
  \renewcommand\thefootnote{}\footnote{#1}%
  \addtocounter{footnote}{-1}%
  \endgroup
}

\input{ex_shared}

\ifpdf
\hypersetup{
  pdftitle={A Herglotz-based integrator for nonholonomic
mechanical systems},
  pdfauthor={Maciel, E., Ortiz, I., Schaerer, C.E.}
}
\fi

\makeatletter
\providecommand*{\diff}%
    {\@ifnextchar^{\DIfF}{\DIfF^{}}}
\def\DIfF^#1{%
    \mathop{\mathrm{\mathstrut d}}%
        \nolimits^{#1}\gobblespace}
\def\gobblespace{%
    \futurelet\diffarg\opspace}
\def\opspace{%
    \let\DiffSpace\!%
    \ifx\diffarg(%
        \let\DiffSpace\relax
    \else
        \ifx\diffarg[%
            \let\DiffSpace\relax
        \else
            \ifx\diffarg\{%
                \let\DiffSpace\relax
            \fi\fi\fi\DiffSpace}

\begin{document}

\maketitle

\footnotesize
\begin{center}
    Polytechnic School - National University of Asuncion\\
San Lorenzo SL 2160 - Paraguay
\end{center} 
\normalsize

\begin{abstract}
We propose a numerical scheme for the time-integration of nonholonomic mechanical systems, both conservative and nonconservative. The scheme is obtained by simultaneously discretizing the constraint equations and the Herglotz variational principle. We validate the method using numerical simulations and contrast them against the results of standard methods from the literature.
\end{abstract}




\section{Introduction}

Numerical integration of differential equations is a very active research area, both from a theoretical perspective as well as a more practical or computational standpoint. This last perspective has had a boost with the advent of digital computers, and we now have a plethora of methods, some of general purposes and some of more ad-hoc nature.

In the last category, we find the {\it structure-preserving integrators}. These are integrators designed to take into account specific properties of the dynamical systems they are meant to solve. Relevant examples of structure-preserving integrators are the symplectic Euler and the symplectic Runge-Kutta, both of which are compatible with the symplectic structure naturally associated with Hamiltonian systems. When the relevant structure is geometric in nature, the structure-preserving integrators are usually known as {\it geometric integrators}. For more details of these and other developments, see \cite{hairer2006geometric}.  

In the last two decades, the differential geometric description of dynamical systems, known as geometric mechanics, and the variational formulation of mechanics, have provided a successful approach for developing geometric integrators. 
Arguably, the seminal result in this context was the fact that, for Hamiltonian systems, discretizing the Hamiltonian principle instead of the Euler-Lagrange equation of motion automatically yields an integrator that preserves the symplectic structure naturally associated with the system, i.e., it yields a {\it symplectic integrator}. This result can be traced back at least to \cite{wendlandt1997mechanical}, and ever since, integrators coming from the discretization of some variational principle are known as {\it variational integrators}. 

The results in \cite{wendlandt1997mechanical} was expanded in \cite{marsden2001discrete} to consider nonconservative systems by discretizing the Lagrange-d'Alembert principle. There, the authors also showed that any symplectic integrator is a variational integrator. Since then, in the following two decades, there has been several generalizations in an attempt to incorporate systems of more general nature, as for instance: coupled multi-body systems, field theories, electrical systems \cite{ding2014higher, leyendecker2008variational, stern2015geometric,  ober2013variational, Galerkin-Variational-2017, MPSW2001Variational-Multisymplectic}; systems with more specific features (as symmetries) \cite{lee2007lie, shen2017lie, demoures2015discrete, marsden1997mechanical, LLM2007Lie-Variational}; interactions with control theory \cite{colombo2012variational, leok2007overview, campos2015high, hussein2006discrete, de2007discrete}; higher-order and asynchronous techniques \cite{LMOW2003AVI, HL2015Spectral-Variational, LS2012Gen-Technique}. For more recent developments, see \cite{Variational-Disipative-2020, variational-2019, Adaptive-Variational-2021, Variational-2020,tran2022multisymplectic, leok2019variational, colombo2022variational, colombo2018variational}. For surveys and reviews, see \cite{Lew2003anoverview, LMO2004Variational-Time-Integrators, Lew2016} and the references therein.

In this work, we construct a numerical integrator for nonconservative and nonholonomically constrained systems. To this end, we combine a criterion for the discretization of nonholonomic constraints (introduced in \cite{cortes2001non}, in the context of Lagrange-d'Alembert principle) and a recently developed integrator for dissipative systems (introduced in \cite{vermeeren2019contact}) coming from a discrete version of the Herglotz's variational principle.  Then we validate our integrator by running simulations of two archetypal mechanical systems in order to contrast it with the outcome of more standard procedures available in the literature.

The organization of this work is as follows. In Section \ref{sec:Background} we give a brief background on geometric mechanics and variational integrators, mainly to fix notation and terminology. In Section \ref{sec:Numerical-Scheme} we provide the technical details of our methodology to get our integrator. In Section \ref{sec:Numerical-experiments} we show the outcome of several numerical simulations, contrasted with a couple of more standard procedures. Finally, in Section \ref{sec:Conclusions} we offer some concluding remarks. We add an appendix containing the explicit equations for the simulations.

\section{Geometric mechanics and Variational integrators}\label{sec:Background}

Geometric mechanics is, roughly speaking, the description of mechanics in the language of differential geometry \cite{holm2011geometric-I, holm2011geometric-II, holm2009geometric-symmetry, marsden2008foundations, marsden1999introduction}. In many interesting cases, the configuration space of a mechanical system with $n\in \mathbb{N}$ degrees of freedom can be modeled as a $n$-dimensional smooth manifold $Q$, and its state space can be modeled by the tangent bundle $TQ$, of the configuration space, or a subset of it. Kinetic energy is given by a Riemannian metric on $Q$, potential energy is given by a function $V\colon Q\to\mathbb{R}$, forces are given by $1$-forms $F\colon Q\to T^*Q$, and the dynamic is encoded by a Lagrangian function $L\colon TQ \to \mathbb{R}$, which typically involves the kinetic and potential energy \cite{bloch2015nonholonomic,marsden1999introduction}.

Constraints of motion are relations among the configuration variables and its derivatives. If those relations involve only the configuration variables, the constraint is called {\it holonomic}, otherwise it is called {\it nonholonomic}. From a geometric point of view, this dichotomy can be seen as follows: holonomically constrained systems are those for which the configuration space is a smooth manifold $Q$, and the state space is the entire tangent bundle $TQ$, which means that, at each configuration $q\in Q$, the entire tangent space $T_qQ$ is allowed for the velocity of the system. On the other hand, nonholonomically constrained systems are those whose configuration space is a smooth manifold $Q$, but the state space is a subset ${\mathcal D}$ of the tangent bundle $TQ$, which means that at each configuration $q\in Q$, only a subset ${\mathcal D}_q\subset T_qQ$ is allowed for the velocity of the system \cite{bloch2015nonholonomic,cortes2002non,lewis2003lagrangian}. For many interesting cases of nonholonomic systems, the constraint subset ${\mathcal D}$ is an affine, or even linear, subbundle of $TQ$ (i.e.: if $Q$ is the configuration manifold, for each point $q\in Q$ the velocity of the system is constrained to a subspace ${\mathcal D}_q$ (linear or affine) of the tangent space $T_qQ$) \cite{Leon1996geometry-lagrangian, Carinena1993geometry-lagrange-multiplier}. In this work we will focus on this kind of constraints. Readers interested in more general cases, including nonlinearity and time dependency may consult the following references \cite{krupkova1997mechanical, krupkova2009nonholonomic-variational, krupkova2010geom-mech-on-nonholonomic-submanifolds, Leon1997nonlinear-nonholonomic, Ranada1994time-dependent-nonholonomic, Budapest1999}.

\begin{figure}
  \centering
  \includegraphics[scale=1.2]{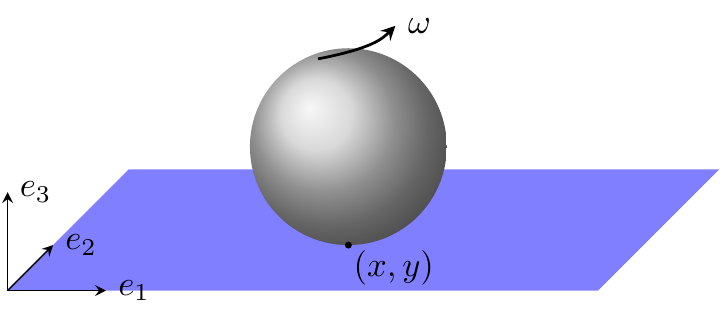}
  \caption{A solid sphere rolling without slipping on a horizontal surface.}
  \label{fig:esfera_que_rueda_sin_deslizarse}
\end{figure}

\begin{example}[\bf Rolling without slipping]
Figure \ref{fig:esfera_que_rueda_sin_deslizarse} illustrates a solid sphere rolling without slipping on a fixed horizontal surface. The configuration of this sphere is fully determined by the point $(x, y)\in \mathbb{R}^2$ of contact between the sphere and the plane, and an element of the special orthogonal group $SO(3)$ (rotations on $\mathbb{R}^3$). Hence, its configuration space is $Q  = \mathbb{R}^2\times SO(3)$, but the nonslipping condition in the contact point imposes a restriction on the velocities; thus, the state space is not the entire tangent bundle $TQ$. This is an archetypal example of a nonholonomically constrained mechanical system.
\end{example}

{\bf Variational formalism:} Given a mechanical system with configuration space $Q$ and constraint subset ${\mathcal D}\subset TQ$, the dynamic (or equation of motion) of the system is modeled by a differential equation on $TQ$ subject to the restriction ${\mathcal D}$. In geometric mechanics, it is common to describe the equation of motion as the extreme of a suitable variational principle. In what follows, we are going to use some standard constructions from the calculus of variations, so for the reader's convenience, let us introduce them for future reference.

Given two points $q_0$ and  $q_f$ in a smooth manifold $Q$, and $t_0<t_f\in \mathbb{R}$, let $\mathcal{C}(t_0, t_f) := \{q\colon[t_0, t_f]\to Q\}$ be the space of smooth paths satisfying $q(t_0) = q_0$ and $q(t_f) = q_f$. A {\bf smooth variation} with fixed endpoints of an element $\overline{q}\in \mathcal{C}(t_0, t_f)$ is a smooth map $q\colon (-\epsilon, \epsilon)\times [t_0, t_f]\to Q$ such that, for each $s\in (-\epsilon, \epsilon)\subset \mathbb{R}$, $q_s(t):= q(s, t)$ for all $t\in [t_0, t_f]$ is an element of $\mathcal{C}(t_0, t_f)$, and $q_0 = \overline{q}$. The {\bf variational vector field} associated to $q$ is the vector field along $\overline{q}$ given by
\[\delta q(\overline{q}(t)):= \frac{\partial}{\partial s}\Big|_{s=0}q(s, t), \forall t\in[t_0, t_f].\]

The simplest scenario for a variational formulation of mechanics is that of a conservative system with holonomic constraints. For this kind of system, the equation of motion is given by Hamilton's Variational Principle, which can be stated as follows.
\begin{definition}[{\bf Hamilton's Principle}]\label{
Hamilton-Principle}
  Let $L\colon TQ\to \mathbb{R}$ be a conservative Lagrangian system. Given two points $q_0$ and  $q_f$ in the configuration space $Q$, and $t_0<t_f\in \mathbb{R}$, the path $\overline{q}(t)\in \mathcal{C}(t_0, t_f)$ followed by the system to go from $q_0$ to $q_f$ is the one that extremizes the action
\[\int_{t_0}^{t_f}L(q, \dot{q})dt.\]
\end{definition}

The optimization mentioned in the Hamilton's Principle is in the sense that, for all smooth variation $q(s, t)$ of $\overline{q}$ in $\mathcal{C}(t_0, t_f)$, with $s$ in a sufficiently small interval $(-\epsilon, \epsilon)\subset \mathbb{R}$, the functional variation of the action vanishes, i.e.:
\[\delta\int_{t_0}^{t_f}L(q, \dot{q})dt := \frac{\partial}{\partial s}\bigg|_{s=0}\left(\int_{t_0}^{t_f}L(q(s, t), \dot{q}(s, t))\right)dt = 0.\]

It is well known that Hamilton's principle is not appropriate for systems that are either nonconservative or with nonholonomic constraints \cite{bloch2015nonholonomic, cortes2002non}. The most widely used variational-like principle to obtain the equation of motion for nonconservative and nonholonomic systems is the Lagrange-d'Alembert principle, together with the so-called nonholonomic principle. This situation can be stated as follows. 

\begin{definition}[{\bf L-A Principle}]
Let $L\colon TQ\to \mathbb{R}$ be the Lagrangian of a system, subject to external forces $F^e\colon Q\to T^*Q$, and nonholonomic constraints ${\mathcal D}\subset TQ$. Given two points $q_0$ and $q_f$ in $Q$, and $t_0< t_f\in \mathbb{R}$, the path $\overline{q}(t)$ followed by the system to go from $q_0$ to $q_f$ is the one that solves
\begin{equation}\label{LA-principle}
    \delta\int_{t_0}^{t_f}L(q, \dot{q})dt + \int_{t_0}^{t_f}F^e \delta q = 0,
\end{equation}
subject to
\begin{equation}\label{nonholonomic-principle}
    \begin{split}
        \delta q(\overline{q}(t))& \in {\mathcal D}_{\overline{q}(t)}, \forall t\in[t_0,t_f],\\
        \dot{\overline{q}}(t)&\in {\mathcal D}_{\overline{q}(t)}, \forall t\in[t_0, t_f].
    \end{split}
\end{equation}
\end{definition}

\begin{remark}
Equation ~\eqref{nonholonomic-principle} is known as the nonholonomic principle. It basically says that both the solution as well as the variational vector field must satisfy the restriction. However, notice that the curves $q_s(t)$, given by the variation $q(s, t)$, are not asked to satisfy the restrictions. Imposing restrictions on those curves poses a different problem, namely, one of optimal control, and in general, it does not provide the equation of motion we are interested in. Readers interested in this dichotomy may consult \cite{bloch2015nonholonomic} and the references therein. 
\end{remark}

{There is an alternative way to handle an interesting spectrum of nonconservative systems by using the Herglotz variational principle \cite{ Georgieva2011Herglotz, ZHANG2019691}. As it turns out, Herglotz principle has a close relationship with contact geometry, which in turn has been found to be a natural arena for several kinds of dissipative systems \cite{de2019contact, Bravetti2017contact, Georgieva2011Herglotz, Bravetti2019Contact-thermodynamics, Bravetti2015Contact-symmetry-thermodynamics, Bravetti2017Contact-Hamiltonian}. In the context of Herglotz principle, the nonconservative nature of mechanical systems with configuration space $Q$ are usually incorporated via a function $L\colon TQ\times \mathbb{R}\to\mathbb{R}$. Due to the aforementioned relationship with contact geometry, we will refer to this function as a {\it contact-type Lagrangian}, to distinguish it from the usual Lagrangian.} 
{
	\begin{remark}
		The term {\it action-dependent Lagrangian} is also used in the literature for Lagrangians in the context of Herglotz principle (see \cite{lazo2018action, de2020review, lopez2022nonsmooth} and the references therein).
\end{remark}}

\begin{definition}[{\bf Herglotz Principle}]\label{Herglotz-Principle}
Let $Q$ be the configuration manifold of a mechanical system whose dynamic is encoded by a {{\it contact-type Lagrangian}} $L\colon TQ\times \mathbb{R}\to \mathbb{R}$ as follows: given $q_0, q_f\in Q$ and $t_0< t_f\in \mathbb{R}$, let ${\mathcal Z}$ be the map that assigns to each $\xi\in {\mathcal C}(t_0, t_f)$, the solution of the initial value problem
\begin{equation}\label{Herglotz-principle}
    \dot{z} = L(\xi, \dot{\xi}, z(t));\quad z(t_0) = z_0.
\end{equation}

\noindent Then, the Herglotz principle establishes that the path $\overline{q}$ followed by the system to go from $q_0$ to $q_f$ is the one that extremizes the functional ${\mathcal Z}(\xi)(t_f)$.
\end{definition}

Table \ref{tab:types_of_systems_and_variational_principles} summarizes some kinds of mechanical systems and the variational principles to obtain their corresponding equations of motion.

\begin{table}[ht]
\centering
\begin{tabular}{l|c|c}
                                                    & \multicolumn{1}{l|}{\textbf{Holonomic}} & \multicolumn{1}{l}{\textbf{Nonholonomic}} \\
\hline
\rotatebox[origin=c]{90}{\hspace{1ex}\textbf{Conservative}\hspace{1ex}}    & \begin{tabular}[c]{@{}c@{}}\textbullet Hamilton's principle\\ $\delta \int_a^b L(q, \dot{q}) \diff t = 0$.\end{tabular}
    & \multicolumn{1}{c|}{\begin{tabular}[c]{@{}c@{}}\textbullet Lagrange-d'Alembert principle\\ $\delta \int_a^b L(q, \dot{q}) \diff t = 0$,\\ $\sum_{k=1}^{n} a_k^j \delta q^k = 0$.\end{tabular}} \\
\hline
\rotatebox[origin=c]{90}{\hspace{1ex}\textbf{Nonconservative}\hspace{1ex}} &
  \begin{tabular}[c]{@{}c@{}}
    \textbullet Lagrange-d'Alembert principle\\ $\delta \int_a^b L(q, \dot{q}) \diff t + \int_a^b F \delta q \diff t = 0$.\\ \\
    \textbullet Herglotz principle\\ $\delta z(b) = 0$. 
  \end{tabular}
    & \multicolumn{1}{c|}{\begin{tabular}[c]{@{}c@{}}
    \textbullet Lagrange-d'Alembert principle\\ $\delta \int_a^b L(q, \dot{q}) \diff t + \int_a^b F^e \delta q \diff t = 0$,\\
    $\sum_{k=1}^{n} a_k^j \delta q^k = 0$.\\ \\
    \textbullet Herglotz principle with constraints\\ $\delta z(b) = 0$,\\
    $\sum_{k=1}^{n} a_k^j \delta q^k = 0$.
    \end{tabular}} \\
\cline{2-3}
\end{tabular}
\caption{Types of mechanical systems and the corresponding variational principles from which to obtain their equations of motion.}
\label{tab:types_of_systems_and_variational_principles}
\end{table}

The variable $z$ in the contact-type Lagrangian for the Herglotz principle allows to incorporate dissipation directly. This makes the Herglotz principle very attractive to obtain integrators that are well suited for dissipative systems, as argued in \cite{vermeeren2019contact}. Next we review the main ideas behind variational integrators.

{\bf Variational integrators:} The main paradigm to get integrators from a variational formulation was introduced in \cite{wendlandt1997mechanical,marsden2001discrete}. There, the authors essentially propose the following procedure: Given a holonomic system $L\colon Q\to \mathbb{R}$, subject to external, nonconservative forces $F^e\colon TQ\to T^*Q$,  discretize the force as $F_d^+, F_d^-\colon Q\times Q\to T^*Q$, as well as the action functional $\int_{t_0}^{t_f}L(q, \dot{q})dt$, namely, let $\{t_i\}_{i=0}^N$ be a partition of the interval $[t_0, t_f]$, and consider 
\[L_d(q_j, q_{j+1})\approx \int_{t_j}^{t_{j+1}}L(q, \dot{q})dt,\]
where $L_d(q_j, q_{j+1})$ is the approximation  given by some quadrature of the integral, and $\{q_j\}_{j=0}^N$ is a sequence of configuration points which approximate the true trajectory $q(t)$ at the times $t_j$, for $j=0,1, 2,\ldots, N$. The sequence $\{q_j\}_{j=0}^N$ is chosen as the one that satisfies the expression
\[\delta\left(\sum_{j=0}^{N-1}L_d(q_j, q_{j+1})\right)+\sum_{j=0}^{N-1}[F_d^-(q_j, q_{j+1})\cdot \delta q_j+F_d^+(q_j, q_{j+1})\cdot \delta q_{j+1}] = 0\]
for all variations $\{\delta q_j\}_{j=0}^N$ vanishing at the endpoints.

For holonomic systems, the variations can be taken in arbitrary directions in each tangent space $T_{q_j}Q$, and thus we get the so-called {\it forced discrete Euler-Lagrange equations:}
\[D_1L_d(q_j, q_{j+1}) + D_2L_d(q_{j-1}, q_j) + F_d^+(q_{j-1}, q_j)+F_d^-(q_j, q_{j+1})= 0,\]
where $D_i$, for $i=1, 2$ denotes the derivative with respect to the $i$-th argument. Under some mild regularity conditions, these equations yield an integrator \[(q_{j-1}, q_{j})\mapsto (q_{j}, q_{j+1}).\]

One of the main results in \cite{marsden2001discrete} is that, for conservative, holonomic systems, this integrator preserves the symplectic structure naturally associated with the given Lagrangian system, i.e., it yields a {\it symplectic integrator}. The symplecticity of the integrator results in very good energy conservation behavior for long-time simulations. 

In this article, we are interested in the approach given in \cite{cortes2001non} to handle nonholonomic systems. There, the authors construct an integrator for nonholonomic systems from the Lagrange-d'Alembert principle, where a main result is a criterion for discretizing the nonholonomic constraint in a compatible way with the discretization of the Lagrangian. To follow up, let us observe that the discrete Lagrangian $L_d$ can be interpreted as a map $L_d \colon Q \times Q \to \mathbb{R}$, which is obtained by means of a discretization map $\Psi \colon Q\times Q \to TQ$ as $L_d = L \circ \Psi$. Then, if ${\mathcal D}\subset TQ$ is the restriction distribution, they propose to take $D:= \Psi({\mathcal D})\subset Q\times Q$ as the discretized restriction. Thus, they get the following integrator:


\begin{equation}\label{Discrete-Const-Forced-EL}
\begin{split}
    D_1L_d(q_j, q_{j+1})+D_2L_d(q_{j-1}, q_j) + F_d^+(q_{j-1}, q_j)+F_d^-(q_j, q_{j+1})& = \lambda_a \Phi^a, \\
    \Phi_d^a(q_j, q_{j+1}) &= 0,
\end{split}
\end{equation}

\noindent where $\lambda_a, a\in\{1, 2,..., m\}$ are Lagrange multipliers and $\Phi^a$ are linearly independent 1-forms defining the annihilator ${\mathcal D}^0\subset T^*Q$ of ${\mathcal D}$.

In \cite{vermeeren2019contact}, the authors derive an integrator for dissipative systems with holonomic constraints by discretizing the Herglotz principle and show that it is compatible with the {\it contact structure} naturally associated to dissipative systems. This compatibility is completely parallel to the one found between symplectic integrators and the symplectic structure of Hamiltonian systems, and thus they named it a {\it contact integrator}.

\section{Construction of the integrator}\label{sec:Numerical-Scheme}
In this section, we construct an integrator for nonholonomically constrained dissipative system following  ideas from \cite{de2019contact}, \cite{vermeeren2019contact}, and \cite{cortes2001non}.

\subsection{Herglotz principle with nonholonomic constraints} Consider a dissipative mechanical system {of contact-type} given by $L\colon TQ\times \mathbb{R}\to \mathbb{R}$ (in the setting of Herglotz principle as in Definition \ref{Herglotz-Principle}), subject to a nonholonomic constraint ${\mathcal D}\subset TQ$.  
Hence, the Herglotz principle is stated as follows: {\it the path $\overline{q}\in {\mathcal C}(t_0, t_f)$ followed by the system is the one that extremizes the functional ${\mathcal Z}(\xi)(t_f)$, among variations tangent to ${\mathcal D}$, and satisfies itself the constraints}. 
This means that, at each time $t\in [t_0, t_f]$, the velocity $\dot{\overline{q}}(t)$ of the path $\overline{q}$ as well as the value $\delta q(t)$ of the variational vector field $\delta q$ must belong to ${\mathcal D}_{\overline{q}(t)}$. 

To express this in a more concrete way, let us consider that the restriction ${\mathcal D}$ is (locally) defined by the vanishing of $m$ functions $\Phi^j\colon TQ\to \mathbb{R}$, linear in the velocities. This means that on a local chart $(q, \dot{q})$ the restrictions are given as the kernel of a matrix $A(q)$; namely, locally we have
\[{\mathcal D} = \{(q, \dot{q}); A(q)\dot{q} = 0\}.\]

\noindent Then, it can be shown  that $\overline{q}$ solves the {\it nonholonomic Herglotz principle} if, and only if, it satisfies the following equations (see \cite{de2019contact} for details): 

\begin{equation}
    \begin{split}
      \frac{d}{dt}\frac{\partial L}{\partial\dot{q}_i}-\frac{\partial L}{\partial q_i} -\frac{\partial L}{\partial\dot{q}_i}\frac{\partial L}{\partial z}  &= \sum_{j=1}^m\lambda_j A_i^j(q),\\
      A(q)\dot{q} &= 0,
    \end{split}
\end{equation}
\noindent for some Lagrange multipliers $\lambda_j$.

 \subsection{Discrete Herglotz Variational Principle}
Consider the discretized Lagrangian $L_d\colon Q^2\times \mathbb{R}^2\to\mathbb{R}$, and the discrete Herglotz principle, stated as follows: given a discrete curve $q = (q_0, q_1, \ldots, q_N)\in Q^{N+1}$, define $z = (z_0, \ldots z_N)\in \mathbb{R}^n$ by $z_0 = 0$ and
  \begin{equation}\label{eq:discrete-action}
  z_{j+1} - z_j = h L_d(q_j, q_{j+1}, z_j, z_{j+1}).
  \end{equation}
  
  Then, a discrete curve $\overline{q} = (\overline{q}_0, \ldots, \overline{q}_N)$ is called a solution of the discrete Herglotz principle if \[\frac{\partial z_{j+1}}{\partial \overline{q_j}}= 0,\]
   for all $j\in \{1, 2, \ldots, N\}$.

It can be shown (see Theorem 1 in \cite{vermeeren2019contact}), that a discrete curve $\overline{q} = (\overline{q}_0, \ldots, \overline{q}_N)$ solves the discrete Herglotz principle if, and only if, it satisfies the discrete generalized Euler-Lagrange equations
\begin{equation} \label{Discrete-Gen-EL}
    \begin{split}
        0 &= D_2L_d(\overline{q}_{j-1}, \overline{q}_j, z_{j-1}, z_j)\frac{1+hD_3L_d(\overline{q}_j, \overline{q}_{j+1}, z_j, z_{j+1})}{1-hD_4L_d(\overline{q}_{j-1}, \overline{q}_j, z_{j-1}, z_j)}\\
        &+ D_1L_d(\overline{q}_j, \overline{q}_{j+1}, z_j, z_{j+1}).
    \end{split}
\end{equation}
Assuming the nondegeneracy condition $D_1 D_2 L_d(q_j, q_{j+1}, z_j, z_{j+1}) \neq 0$, equation~\eqref{Discrete-Gen-EL} yields an integrator $(q_{j-1}, q_j) \mapsto (q_{j}, q_{j+1})$, for the unconstrained system $L\colon TQ\times \mathbb{R}\to \mathbb{R}$. In \cite{vermeeren2019contact}, it is proved that this integrator preserves the natural contact structure associated to the  {discrete contact-type} Lagrangian $L_d$, hence the authors called it a {\it contact integrator}.
 {\begin{remark}
A simpler (and perhaps more natural) discretization of $TQ\times \mathbb{R}$ would be $Q^2\times\mathbb{R}$, as was used in \cite{anahory2021geometry}. This would lead to a discretized action-type Lagrangian $L_d$ without the $z_{j+1}$ argument and without the $D_4L_d$ term in ~\eqref{Discrete-Gen-EL}. As pointed out in \cite{vermeeren2019contact}, for the theoretical development of contact integrators, there is no need for this extra argument in $L_d$. However, the discretization without the $z_{j+1}$ argument is of first order, whereas the one with  $z_{j+1}$ is of second order, and this is relevant for numerical behavior. In our work, for the numerical experiments, we will use both alternatives.
\end{remark}}

\subsection{Discretization of the Constraint}
Now, let us bring into play the constraint ${\mathcal D}\subset TQ$. Let us first recall the method introduced in \cite{cortes2001non} to incorporate a nonholonomic constraint into the integrator derived from the Lagrange-d'Alembert principle, given a mechanical system $L\colon TQ\to \mathbb{R}$. Since the discretization of the Lagrangian yields a map $L_d\colon Q\times Q\to \mathbb{R}$, the discrete constraint must be a submanifold $\Delta_d\subset Q\times Q$, which satisfies $(q, q)\in \Delta_d$ for all $q\in Q$,  {(i.e., $\Delta_d$ has to contain the diagonal of $Q\times Q$).} This submanifold imposes a restriction to the motion in the sense that the discrete curve $q = (q_0, q_1, \ldots, q_N)$ must satisfy that $(q_j, q_{j+1})\in \Delta_d$, and the variation $\delta q$ must satisfy the condition $\delta q(q_j)\in {\mathcal D}_{q_j}$, i.e., the allowed variations are those who respect the continuous constraint ${\mathcal D}\subset TQ$. In order to construct $\Delta_d$ in a compatible way with the discrete Lagrangian, it is assumed that there is a discretizing map $\Psi\colon Q\times Q\to TQ$ such that $L_d = L\circ \Psi$, and then, $\Delta_d:= \{(q, q')\in Q\times Q; \Psi(q, q')\in {\mathcal D}\}$. This means that, if ${\mathcal D}$ is (locally) defined by the vanishing of the functions $\Phi^c\colon TQ\to \mathbb{R}$, then $\Delta_d$ is (locally) defined by the vanishing of the functions $\Phi^c_d:= \Phi^c\circ\Psi$.

In our situation, as we have a  {contact-type} Lagrangian $L\colon TQ\times \mathbb{R}\to \mathbb{R}$, we might consider a discretizing map of the form $\Psi = \Psi_1\otimes \Psi_2\colon Q^2\times \mathbb{R}^2\to TQ\times \mathbb{R}$, such that $\Psi(q, q', z, z') = (\Psi_1(q, q'), \Psi_2(z, z'))$. Thus, we want our discrete restriction to be compatible with the discretizing map $\Psi_1\colon Q^2\to TQ$, in the same way as discussed in the previous paragraph.  

More concretely, suppose as before that on a chart $(q, \dot{q})$ of $TQ$ our restriction ${\mathcal D}$ is given in the form
\[{\mathcal D} = \{(q, \dot{q})\in TQ; A(q)\dot{q} = 0\},\]
for some matrix $A(q)$. Let $(q_d, \dot{q}_d)\in TQ$ be the image of a point $(q, q')\in Q^2$ under the map $\Psi_1$. Then we define the discrete restriction map as $A_d(q, q') = A(q_d)\dot{q}_d$, and the discrete restriction as
\[\Delta_d := \{(q, q')\in Q^2; A_d(q, q') = 0\}.\]

Therefore, the discrete generalized Euler-Lagrange equations with nonholonomic constraints linear in the velocities are
\small
\begin{equation}\label{constrained-discrete-GEL-eq}
\begin{split}
  D_1 L_d(q_j, q_{j+1}, z_j, z_{j+1}) + D_2 L_d(q_{j-1}, q_j, z_{j-1}, z_j) \frac{1 + h D_3 L_d(q_j, q_{j+1}, z_j, z_{j+1})}{1 - h D_4 L_d(q_{j-1}, q_j, z_{j-1}, z_j)} &= \lambda A(q_j), \\
  A_d(q_j, q_{j+1}) &= 0,
\end{split}
\end{equation}
\normalsize
which correspond to the equations of a contact integrator with nonholonomic constraints linear in the velocities.

\section{Numerical simulations and comparison}\label{sec:Numerical-experiments}
In this section, we will apply our integrator to two nonholonomic mechanical systems, namely the Foucault pendulum and the falling rolling disk.

\subsection{Foucault pendulum}
We analyze the Foucault pendulum with Rayleigh dissipation. We derive its Lagrangian, both for the Lagrange-d'Alembert as well as for the Herglotz formulation, in order to compare our integrator with the one coming from the Lagrange-d'Alembert principle.

\begin{figure}
  \centering
  \begin{subfigure}[b]{0.45\textwidth}
  \centering
    \includegraphics{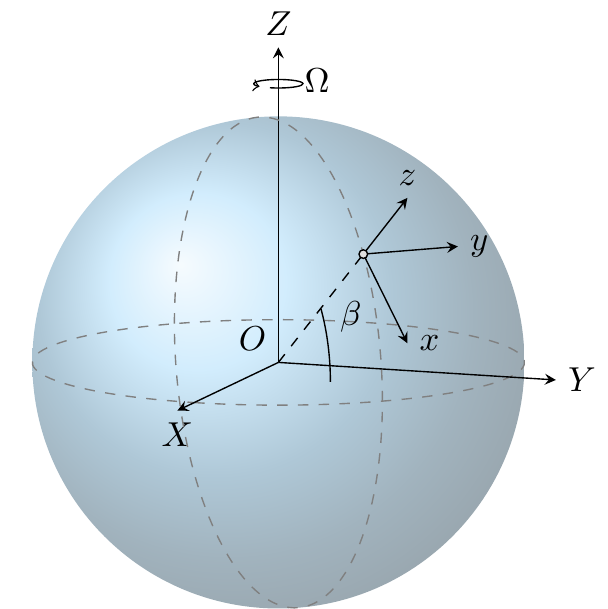}
  \end{subfigure}
  \begin{subfigure}[b]{0.45\textwidth}
  \centering
    \includegraphics{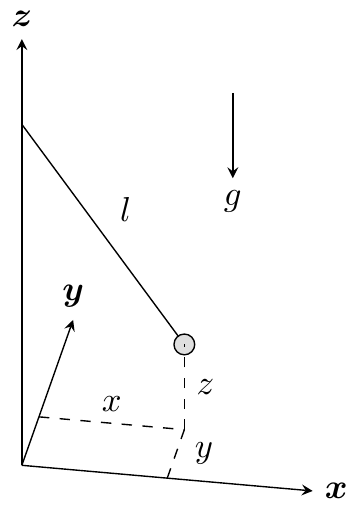}
  \end{subfigure}
  \caption{Foucault pendulum and its inertial and noninertial frames.}
  \label{fig:pendulo_de_foucault}
\end{figure}

\noindent{\bf General setting and simplifications:} The Foucault pendulum consists of a pendulum of length $l$ and mass $m$ located at latitude $\beta$ on Earth's surface. We may consider here a Rayleigh dissipation with parameter $\alpha$. The plane of oscillation does not rotate in a reference $(X, Y, Z)$ fixed in space, hence as Earth rotates, this plane rotates with respect to a reference $(x, y, z)$ attached to Earth (see Figure \ref{fig:pendulo_de_foucault}). 

To model this problem, we consider the inertial frame $(X, Y, Z)$ with its origin at the center of the Earth and $Z$ passing through the north pole. Hence, the angular velocity of Earth, $\Omega$, points along $Z$. On the other hand, the noninertial frame $(x, y, z)$ is such that $x$ points along a meridian in the south direction, $y$ points to the east along the parallel $\beta$ and $z$ coincides with the vertical at the pendulum location. Hence, the vector position $r = (x, y, z)$ of a particle in the noninertial frame satisfies the relation $\omega = \frac{r\times \dot{r}}{\|r\|^2}$, where $\omega = (\Omega \cos \beta, 0, -\Omega \sin \beta)$ is the angular velocity of the particle in the noninertial frame. 

As the pendulum forms a small angle $\phi$ in its oscillatory motion, the
coordinates $x, y$ are of order $l \phi$, whereas $z$ is of order $l \phi^2$
and so it is negligible. Thus, we may consider the
movement of the mass pendulum in the plane ${z=0}$ and take $q = (x, y)$ as generalized coordinates. With these considerations, the kinetic energy $K(\dot{q})$ and potential energy $V(q)$, in terms of the mass $m$, the length $l$ and the gravitational acceleration $g$, are given by:
\begin{equation}\label{energy-Foucault}
  K(\dot{q}) = \frac{1}{2}m(\dot{x}^2+\dot{y}^2)\qquad\text{and}\qquad V(q)= \frac{1}{2} m \frac{g}{l} (x^2 + y^2),
\end{equation}
while the relation $\omega = \frac{r\times \dot{r}}{\|r\|^2}$ reads
\begin{equation}\label{nonhol-restriction}
  - y \dot{x} + x \dot{y} + \Omega \sin \beta (x^2 + y^2) = 0,
\end{equation}
which is a nonholonomic constraint for the system.\\

\noindent{\bf Lagrange-d'Alembert description:} In this context, the Lagrangian is $L(q, \dot{q}) = K(\dot{q})-V(q)$, which according to our previous computations reads
\begin{equation}\label{eq:LD-Foucault}
  L(q, \dot{q}) = \frac{1}{2} m (\dot{x}^2 + \dot{y}^2) - \frac{1}{2} m
  \frac{g}{l} (x^2 + y^2),
\end{equation}
subject to the constraint~\eqref{nonhol-restriction}. The damping is modeled as an external force $F(q, \dot{q})=-\alpha m\dot{q}$.

To derive the integrator for this case, we use equation~\eqref{Discrete-Const-Forced-EL} with a linear-order quadrature for both the Lagrangian and the external force, while the discrete constraint is obtained using the discretizing map $\Psi$ on equation \eqref{nonhol-restriction} (see~\cite{cortes2001non}). The resulting equations are
\begin{equation}
\begin{split}
  \frac{- x_{j+1} + 2 x_j - x_{j-1}}{h} - h \frac{g}{l} x_j - \alpha ( x_{j+1} - x_j ) + \lambda_1 \frac{y_j}{m} &= 0 \\
  \frac{- y_{j+1} + 2 y_j - y_{j-1}}{h} - h \frac{g}{l} y_j - \alpha ( y_{j+1} - y_j ) - \lambda_1 \frac{x_j}{m} &= 0 \\
  - y_j \frac{x_{j+1} - x_j}{h} + x_j \frac{y_{j+1} - y_j}{h} + \Omega \sin \beta (x_j^2 + y_j^2) &= 0,
\end{split}
\end{equation}
where $\lambda_1$ is a Lagrange multiplier.

\vspace{12px}

\noindent{\bf Herglotz description:} Here, the Lagrangian is taken as $L(q, \dot{q}, z) = K(\dot{q}) - V(q) - \alpha z$, which, according to our previous analysis becomes
\begin{equation} \label{eq:H-Foucault}
  L(q, \dot{q}, z) = \frac{m}{2} (\dot{x}^2 + \dot{y}^2) - 
  \frac{m g}{2 l} \left( x^2 + y^2 \right) - \alpha z,
\end{equation}
also subject to the nonholonomic constraint~\eqref{nonhol-restriction}.

To derive the Contact integrator for this system, we use a linear-order approximation in equations \eqref{constrained-discrete-GEL-eq}, i.e.
\begin{equation}
  z_{j+1} - z_j = h L(x_j, x_{j+1}, z_{j}, z_{j+1}),
\end{equation}
and the constraints are discretized in the same way as in \cite{cortes2001non}. The resulting
equations are
\begin{equation}
\begin{split}
  \frac{- x_{j+1} + 2 x_j - x_{j-1}}{h^2} - \frac{g}{l} x_j - \alpha \left( \frac{x_{j} - x_{j-1}}{h} - \frac{h}{2} \frac{g}{l} x_j \right) + \lambda_1 \frac{y_j}{m} &= 0 \\
  \frac{- y_{j+1} + 2 y_j - y_{j-1}}{h^2} - \frac{g}{l} y_j - \alpha \left( \frac{y_{j} - y_{j-1}}{h} - \frac{h}{2} \frac{g}{l} y_j \right) - \lambda_1 \frac{x_j}{m} &= 0 \\
  - y_j \frac{x_{j+1} - x_j}{h} + x_j \frac{y_{j+1} - y_j}{h} + \Omega \sin \beta (x_j^2 + y_j^2) &= 0,
\end{split}
\end{equation}
where $\lambda_1$ is a Lagrange multiplier.

\vspace{12px}

As can be noted, the discretization of the constraints in both integrators LA and Contact are the same.

\subsection{Falling disk} \label{sec:falling_disk}

\begin{figure}
  \centering
  \includegraphics[scale=1.2]{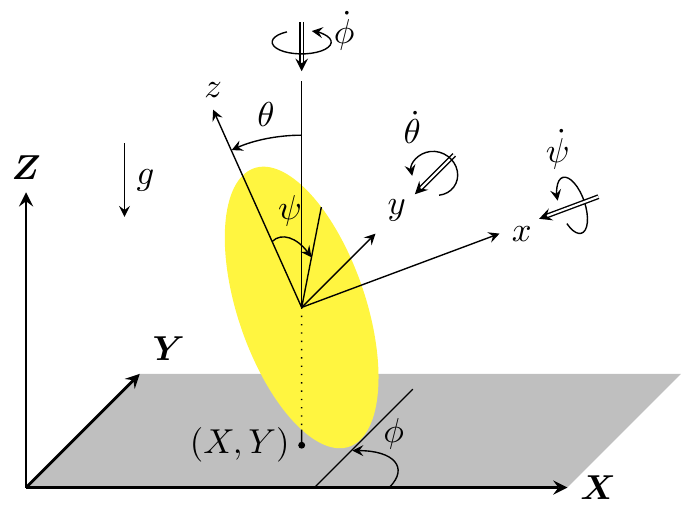}
  \caption{Falling disk and its state variables.}
  \label{fig:disco_rodante_que_cae}
\end{figure}

The rolling falling disk is one of the archetypal examples of nonholonomically constrained mechanical systems, frequently used as a benchmark for numerical methods \cite{bloch2015nonholonomic, modin2017makes}. It consists of an idealized (without thickness) homogeneous disk rolling without slipping over a fixed horizontal plane, subject to gravity (see Figure \ref{fig:disco_rodante_que_cae}). Depending on the research interest, we may also consider external forces acting upon the disk, either dissipative forces or control forces. To derive its Lagrangian, we will follow the formulation given in \cite{paris2002disk}. 

The configuration space of this system can be identified as $Q= \mathbb{R}^2\times SO(3)$. As generalized coordinates, we can take a point $(X, Y)\in \mathbb{R}^2$ describing the orthogonal projection of the center of the disk over the plane $\{Z = 0\}$ (see Figure \ref{fig:disco_rodante_que_cae}), and an element of $SO(3)$ in terms of Eulerian angles $(\theta, \phi, \psi)$, defined as follows: we consider a right-handed noninertial orthogonal reference frame $(x, y, z)$ whose origin is attached to the center of the disk, the $x$ axis being orthogonal to the plane of the disk and the $y$ axis remaining in the plane of the disk and parallel to the plane $\{Z = 0\}$. Then, $\theta$ is the angle between the axes $z$ and $Z$, while $\phi$ is the angle between the axes $y$ and $X$, and $\psi$ is the angle of rotation of the disk with respect to the $x$ axis. 

In the Lagrange-d'Alembert formulation, the Lagrangian for this system is $L = T- V$, where $T$ is the kinetic energy and $V$ is the potential energy. For the kinetic energy, we have the expression
\begin{equation}
  T = \frac{1}{2}m\hat{v}^2 + \frac{1}{2} [I_A \omega_x^2 + I_T (\omega_y^2 + \omega_z^2)],
\end{equation}
where $\hat{v}$ is the linear velocity of the center of the disk, $I_A$ and $I_T$ are the moment of inertia with respect to $x$ axis, and the $y$ and $z$ axes, respectively, while $\omega_i$, for $i\in\{x, y, z\}$ are the angular velocity with respect to the $i$-th axis. In terms of the generalized velocities associated to the Euler angles, the angular velocities are given by
\begin{equation}
\begin{split}
  \omega_x &= -\dot{\psi} + \dot{\phi} \sin\theta, \\
  \omega_y &= -\dot{\theta}, \\
  \omega_z &= \dot{\phi} \cos\theta.
\end{split}
\end{equation}
On the other hand, the potential energy is given by
\[V = m g R \cos \theta,\]
where $m$ is the mass of the disk, $g$ is the gravitational acceleration, and $R$ is the radius of the disk. Hence, the complete Lagrangian is
\begin{equation}\label{eq:Disk-lagrangian-LD-form}
\begin{split}
  L = T(\dot{q}) - V(q) &= \frac{1}{2} m [\dot{X}^2 + \dot{Y}^2 + R^2 \sin^2\theta \dot{\theta}^2] \\
  &+ \frac{1}{2} \left[ I_A \left( \dot{\psi} - \dot{\phi} \sin\theta \right)^2 + 
  I_T \left( \dot{\theta}^2 + \dot{\phi}^2 \cos^2\theta \right) \right] - m g R \cos\theta.
\end{split}
\end{equation}
Finally, the non-slipping condition imposes restrictions on the velocities given by

\begin{equation} \label{eq:disco_que_cae_restricciones}
\begin{split}
  \dot{X} &= - R \cos\theta \sin\phi \dot{\theta} - R \sin\theta \cos\phi \dot{\phi} + R \cos\phi \dot{\psi}, \\
  \dot{Y} &= R \cos\theta \cos\phi \dot{\theta} - R \sin\theta \sin\phi \dot{\phi} + R \sin\phi \dot{\psi}.
\end{split}
\end{equation}
Let us consider that the motion of the disk undergoes a Rayleigh dissipation and also that there is some external nonconservative force applied upon the disk. The Herglotz description allows us to incorporate these features into the Lagrangian of the system by modifying the expression given in ~\eqref{eq:Disk-lagrangian-LD-form} in the following way
\begin{equation}
  L_H(t, q, \dot{q}, z) = L -\alpha z + F(t)q,
\end{equation}
where $L$ is the Lagrangian given by ~\eqref{eq:Disk-lagrangian-LD-form}, $\alpha$ is the dissipation parameter, and $F(t)$ is the nonconservative external force. The restriction of the system, given by ~\eqref{eq:disco_que_cae_restricciones} are unchanged in the Herglotz formulation.

\subsection{Numerical Experiments: Foucault Pendulum} \label{sec:numerical_experiments}
For the Foucault pendulum, we perform two simulations, using as reference the $4^{th}$ order Runge-Kutta-Fehlberg method, to compare the behavior of both integrators. In both simulations, we use the configuration of the original experiment held in the Observatory of Paris, in 1851, namely, we consider a pendulum of mass $m = 28$ kg, length $l = 67$ m, and latitude $\beta = 49^\circ$. We consider two different values of the dissipation parameter $\alpha$. 

\begin{figure}
  \centering
  \begin{subfigure}[b]{0.65\textwidth}
    \centering
    \includegraphics[width=\textwidth]{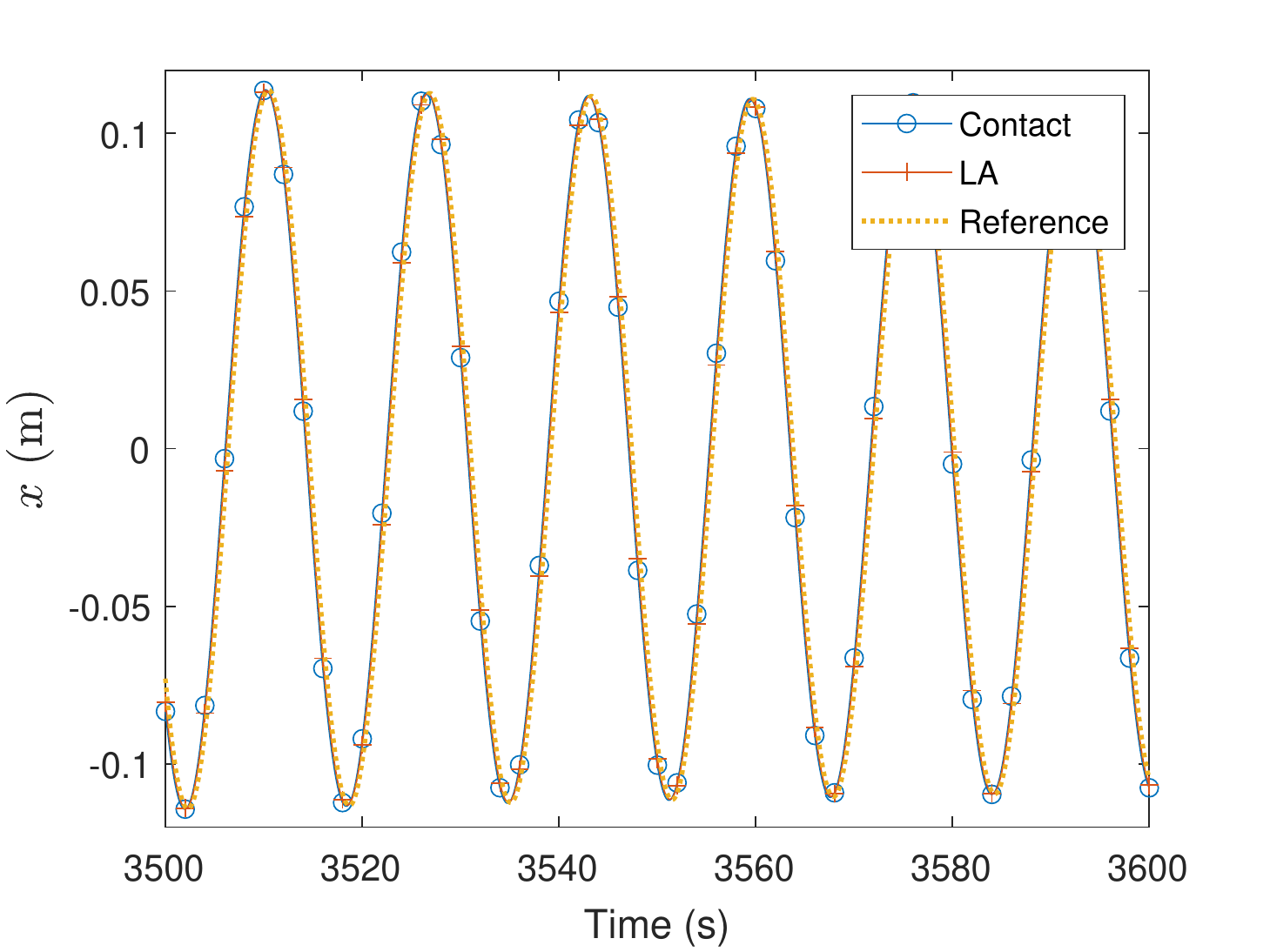}
  \end{subfigure} \\
  \begin{subfigure}[b]{0.65\textwidth}
    \centering
    \includegraphics[width=\textwidth]{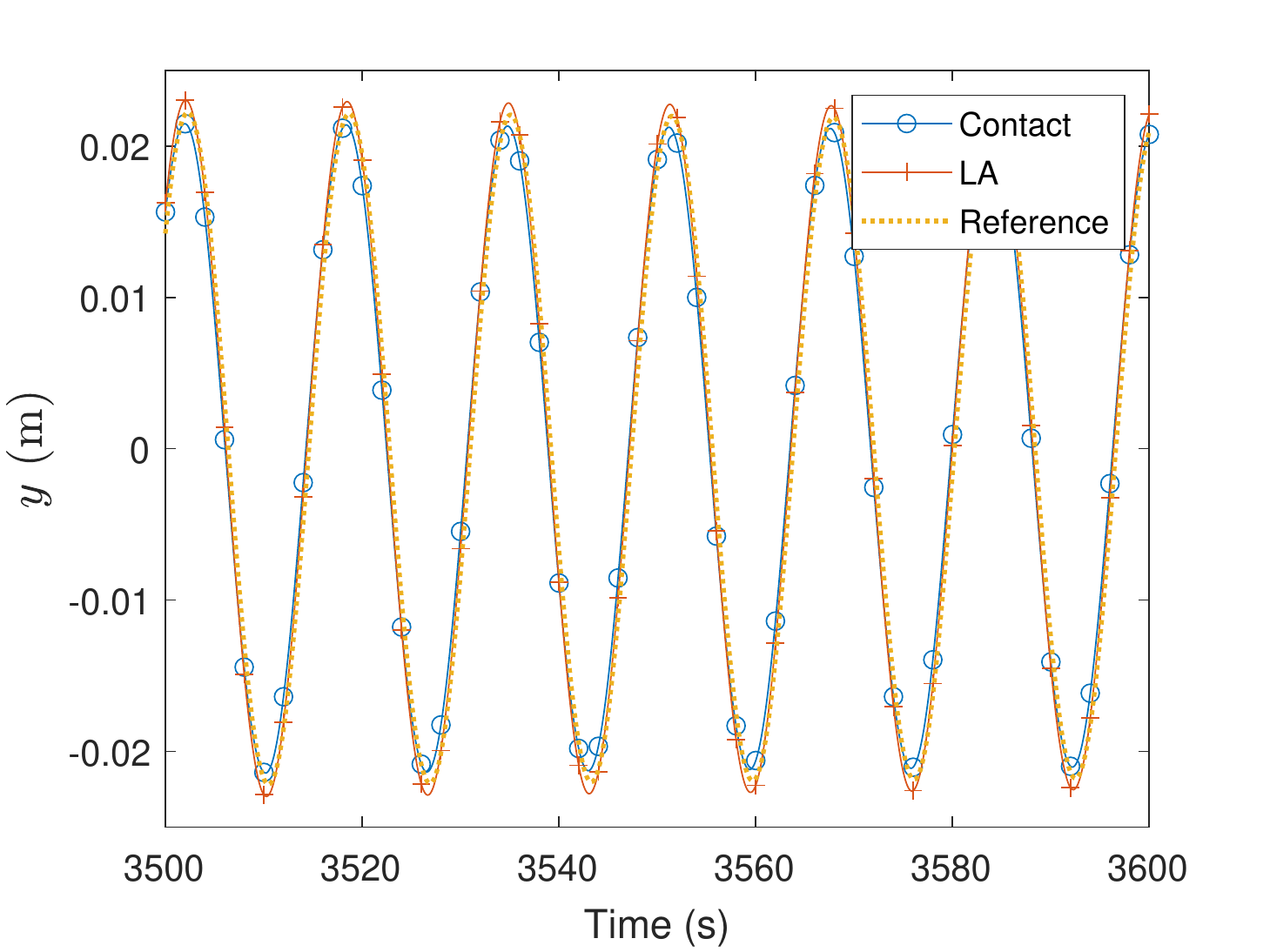}
  \end{subfigure}
  \caption{Simulation 1 for the Foucault Pendulum with initial conditions $q(0) = (0,l/100)$, $\dot{q}(0) = (0,0)$. Dissipation parameter $\alpha = 1 \times 10^{-3}$. The graphics show the values of the $x$ and $y$ coordinates as functions of time for the last $100$ seconds of simulation.}
  \label{fig:foucault_sim_1_trayectories}
\end{figure}

\begin{figure}
  \centering
  \begin{subfigure}[b]{0.65\textwidth}
    \centering
    \includegraphics[width=\textwidth]{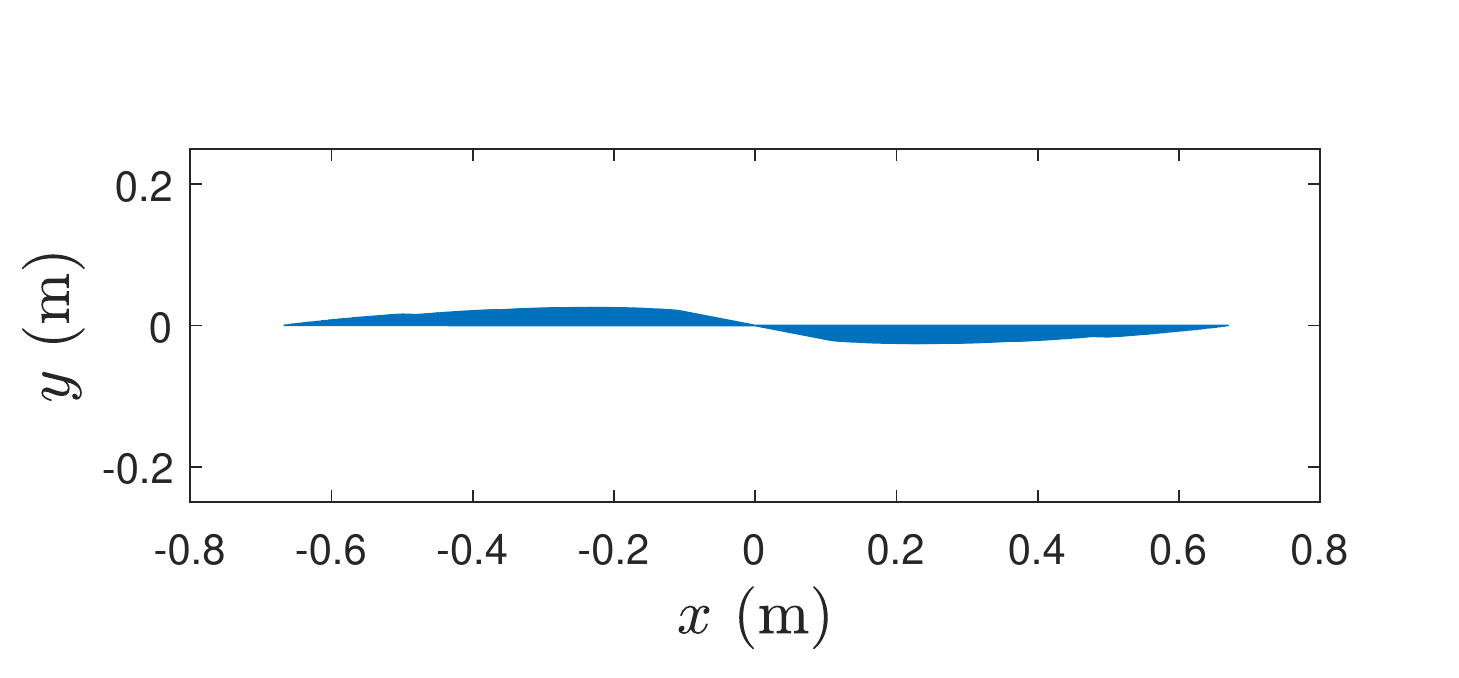}
  \end{subfigure} \\
  \begin{subfigure}[b]{0.65\textwidth}
    \centering
    \includegraphics[width=\textwidth]{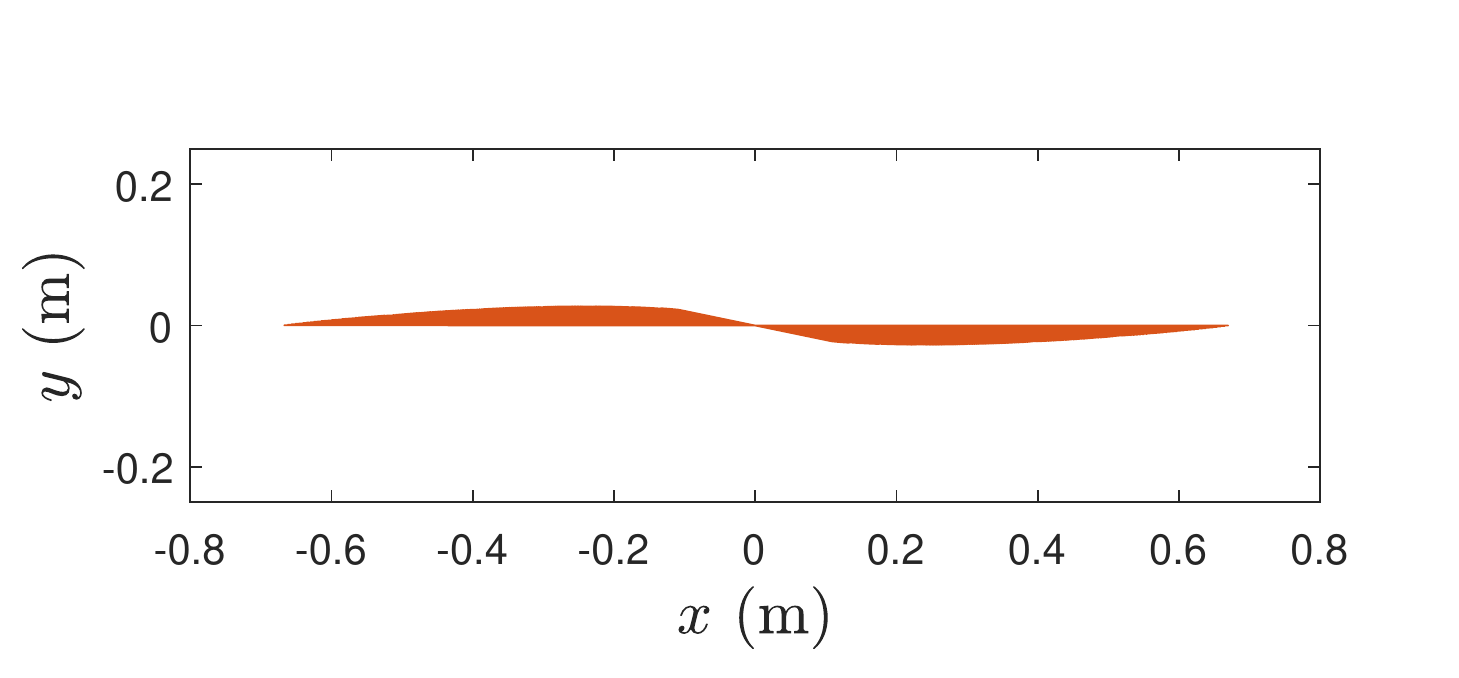}
  \end{subfigure}
  \begin{subfigure}[b]{0.65\textwidth}
    \centering
    \includegraphics[width=\textwidth]{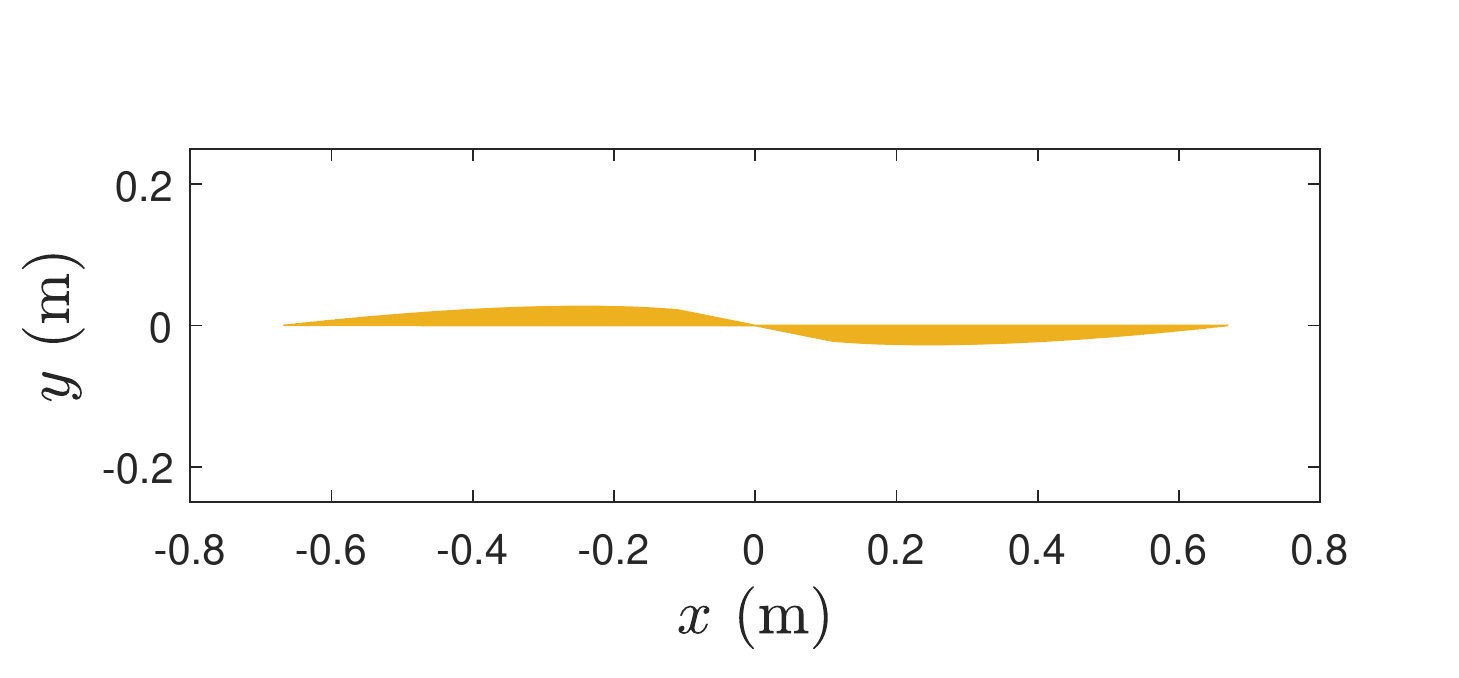}
  \end{subfigure}
  \caption{Simulation 1 for the Foucault Pendulum with initial conditions $q(0) = (0,l/100)$, $\dot{q}(0) = (0,0)$. Dissipation parameter $\alpha = 1\times 10^{-3}$. The graphics show the projection over the horizontal plane of the oscillation plane in the interval $0\leq t\leq 3600$ seconds.}
  \label{fig:foucault_sim_1_oscillation_plane}
\end{figure}

\begin{figure}
  \centering
  \begin{subfigure}[b]{0.65\textwidth}
    \centering
    \includegraphics[width=\textwidth]{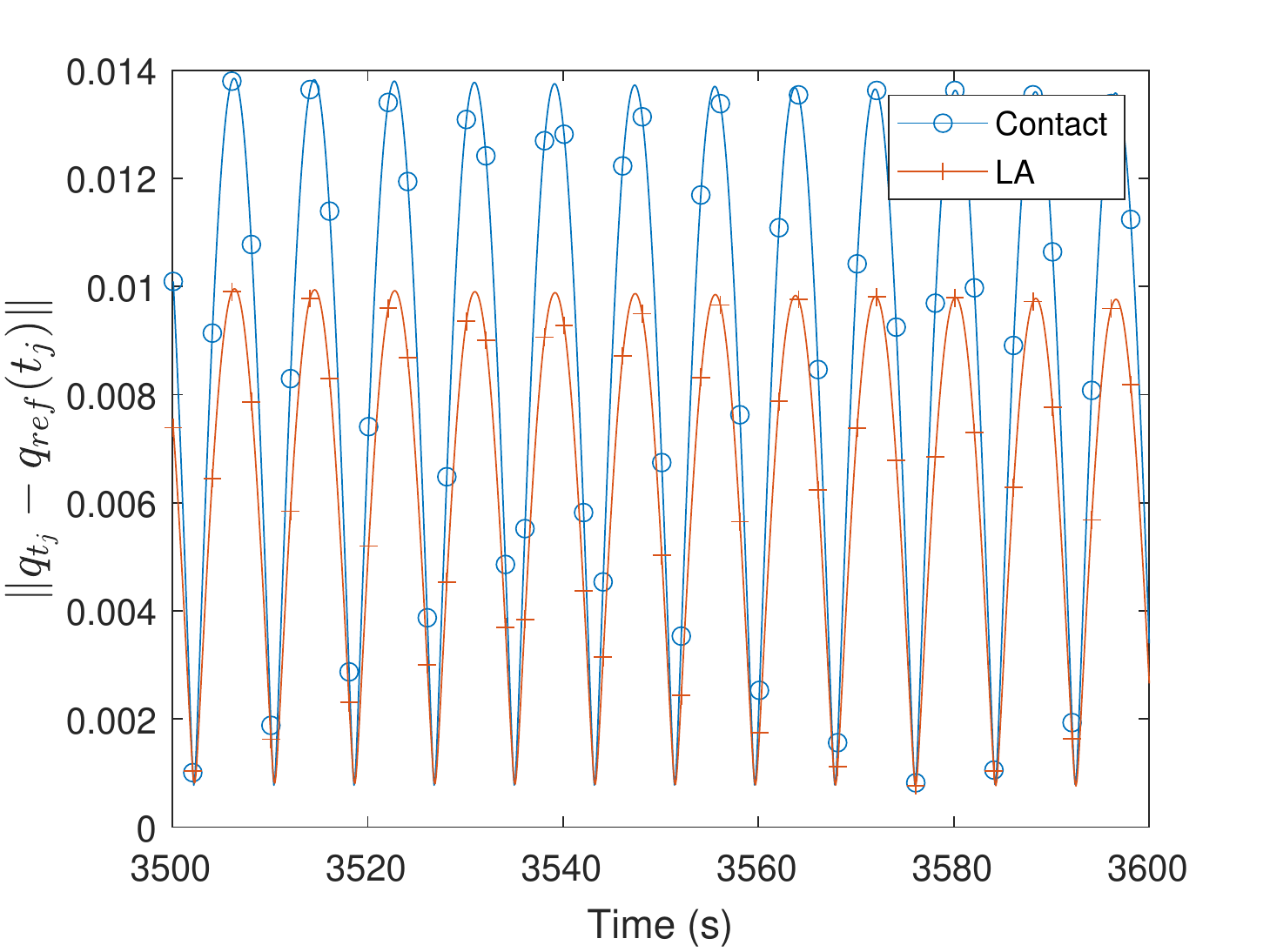}
    \caption{Trajectory error}
  \end{subfigure} \\
  \begin{subfigure}[b]{0.65\textwidth}
    \centering
    \includegraphics[width=\textwidth]{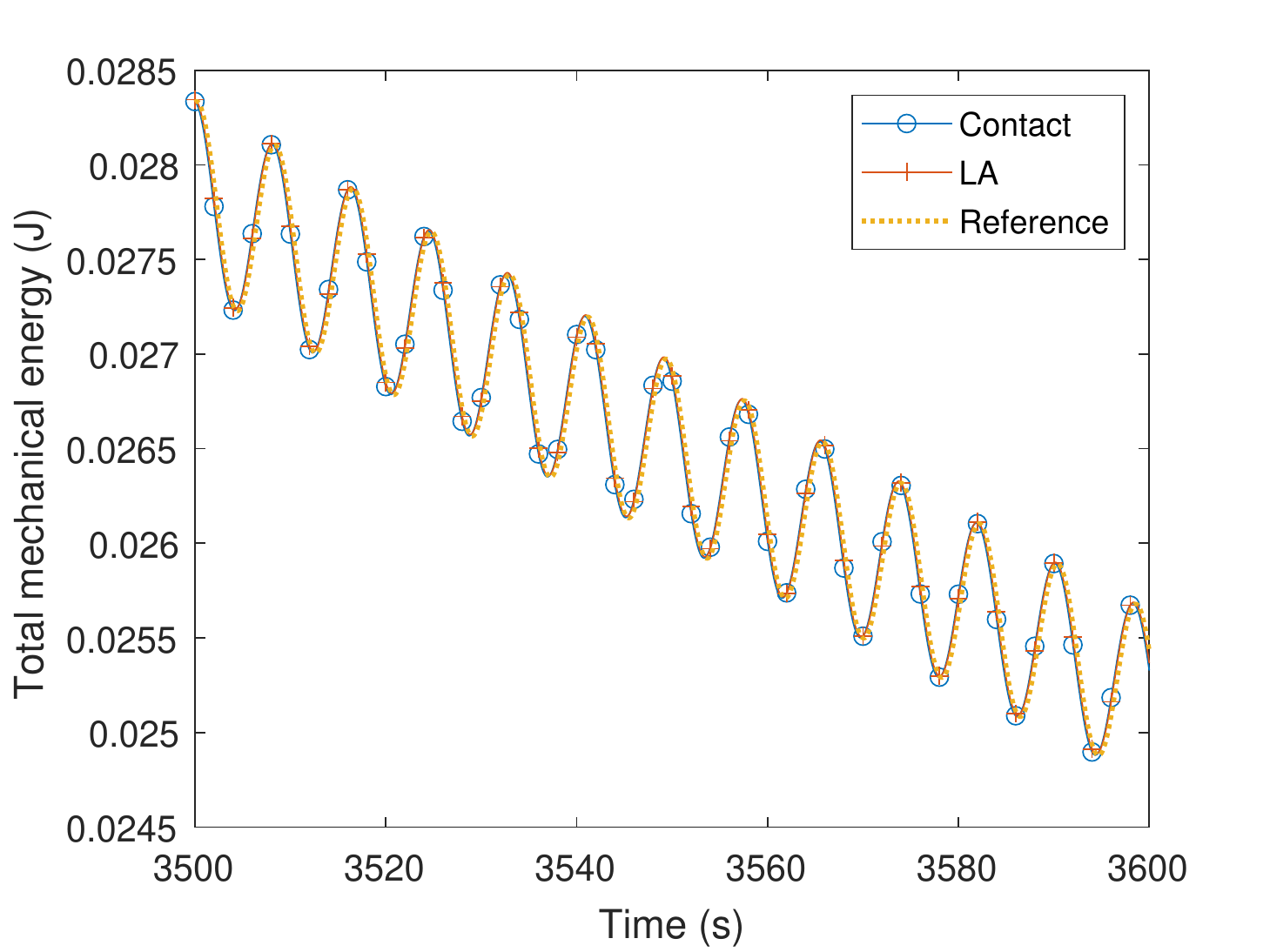}
    \caption{Energy}
  \end{subfigure}
  \caption{Simulation 1 for the Foucault Pendulum with initial conditions $q(0) = (0,l/100)$, $\dot{q}(0) = (0,0)$. Dissipation parameter $\alpha = 1\times 10^{-3}$.  The graphics show the errors of the trajectory and the energy for the last $100$ seconds of simulation.}
  \label{fig:foucault_sim_1_err_energy}
\end{figure}

\vskip 2mm
\noindent{\bf Experiment 1:} Here we consider $\alpha = 1\times 10^{-3}$ for both the Contact and the LA integrators. For this setting, both integrators display indistinguishable behavior in terms of the trajectory (see Figure \ref{fig:foucault_sim_1_trayectories} and Figure \ref{fig:foucault_sim_1_oscillation_plane}) as well as in energy dissipation, as can be seen in part b) of Figure \ref{fig:foucault_sim_1_err_energy}.  However, we observe that the LA integrator do have a better outcome in terms of the error $\|q_{tj} - q_{ref}(t_j)\|_2$, where $q_{t_j}$ is the approximation given by the integrator and $q_{ref}(t_j)$ is the reference solution, both at time $t=t_j$. This can be seen in the part a) of Figure \ref{fig:foucault_sim_1_err_energy}.

\begin{figure}
  \centering
  \begin{subfigure}[b]{0.65\textwidth}
    \centering
    \includegraphics[width=\textwidth]{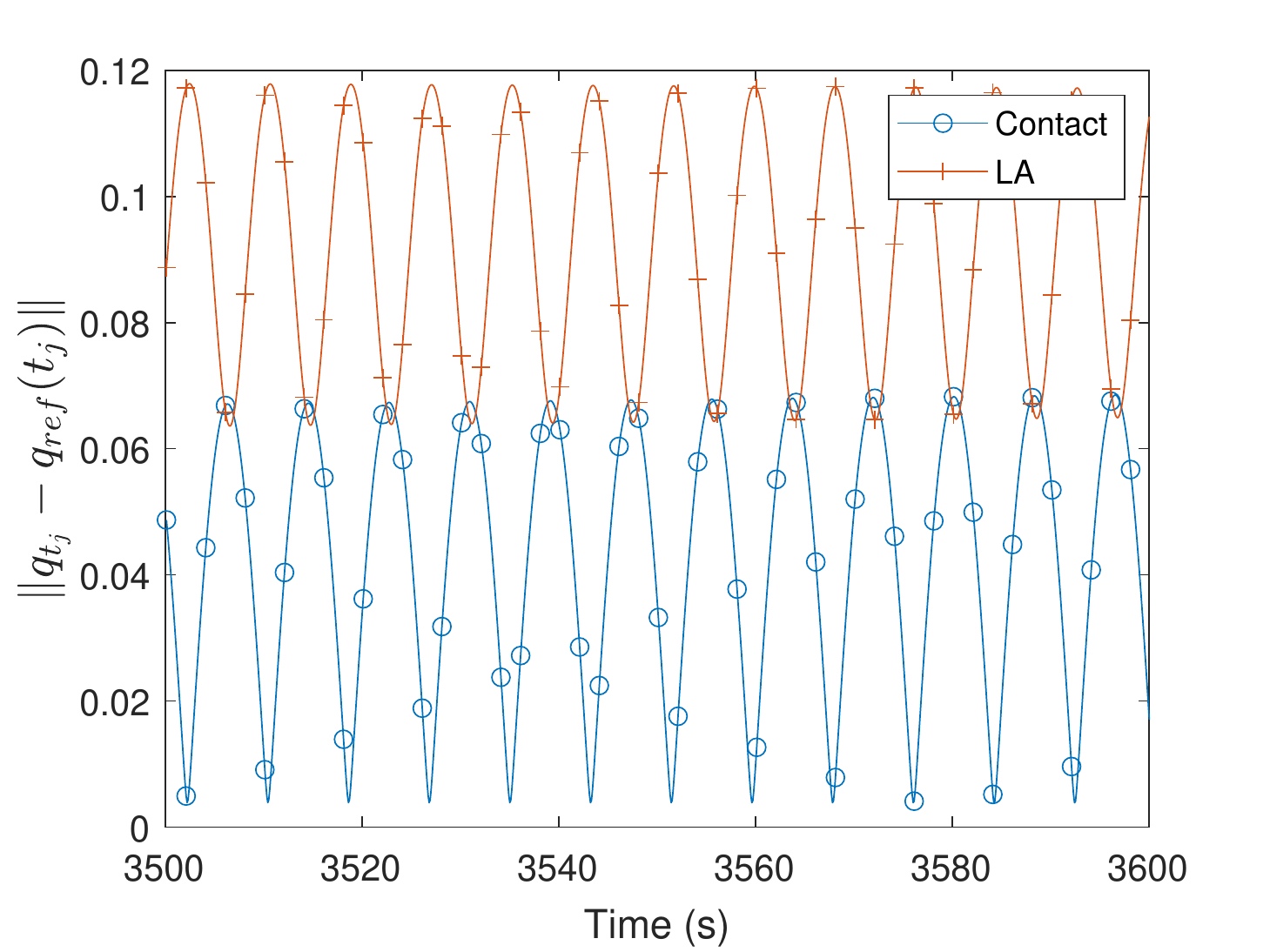}
  \end{subfigure} \\
  \begin{subfigure}[b]{0.65\textwidth}
    \centering
    \includegraphics[width=\textwidth]{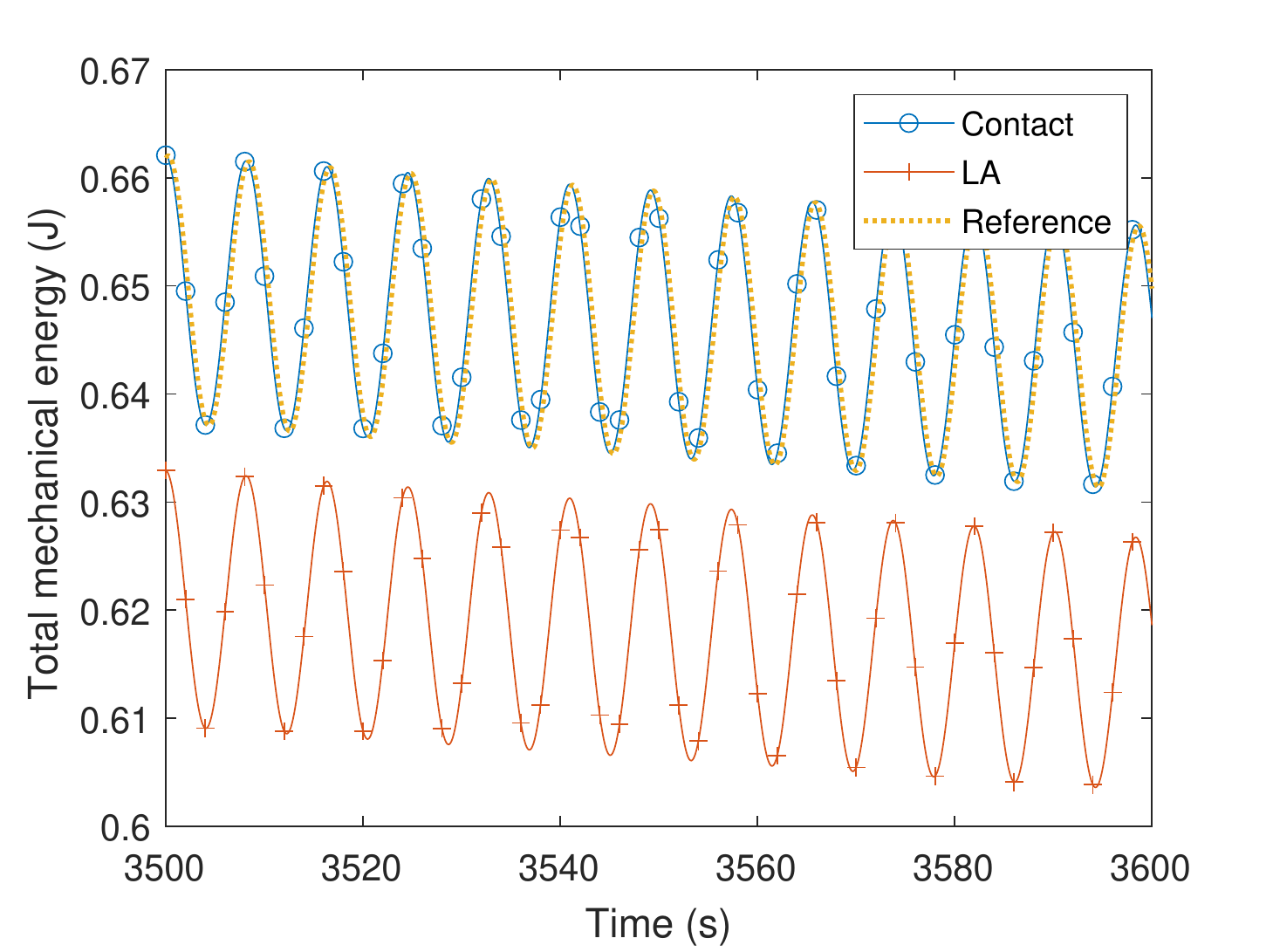}
  \end{subfigure}
  \caption{Simulation 2 for the Foucault Pendulum with initial conditions $q(0) = (0,l/100)$, $\dot{q}(0) = (0,0)$. Dissipation parameter $\alpha = 1\times 10^{-4}$. The graphics show the errors of the trajectory and the energy for the last $100$ seconds of simulation. }
  \label{fig:foucault_sim_2_err_energy}
\end{figure}

\begin{figure}
  \centering
  \begin{subfigure}[b]{0.65\textwidth}
    \centering
    \includegraphics[width=\textwidth]{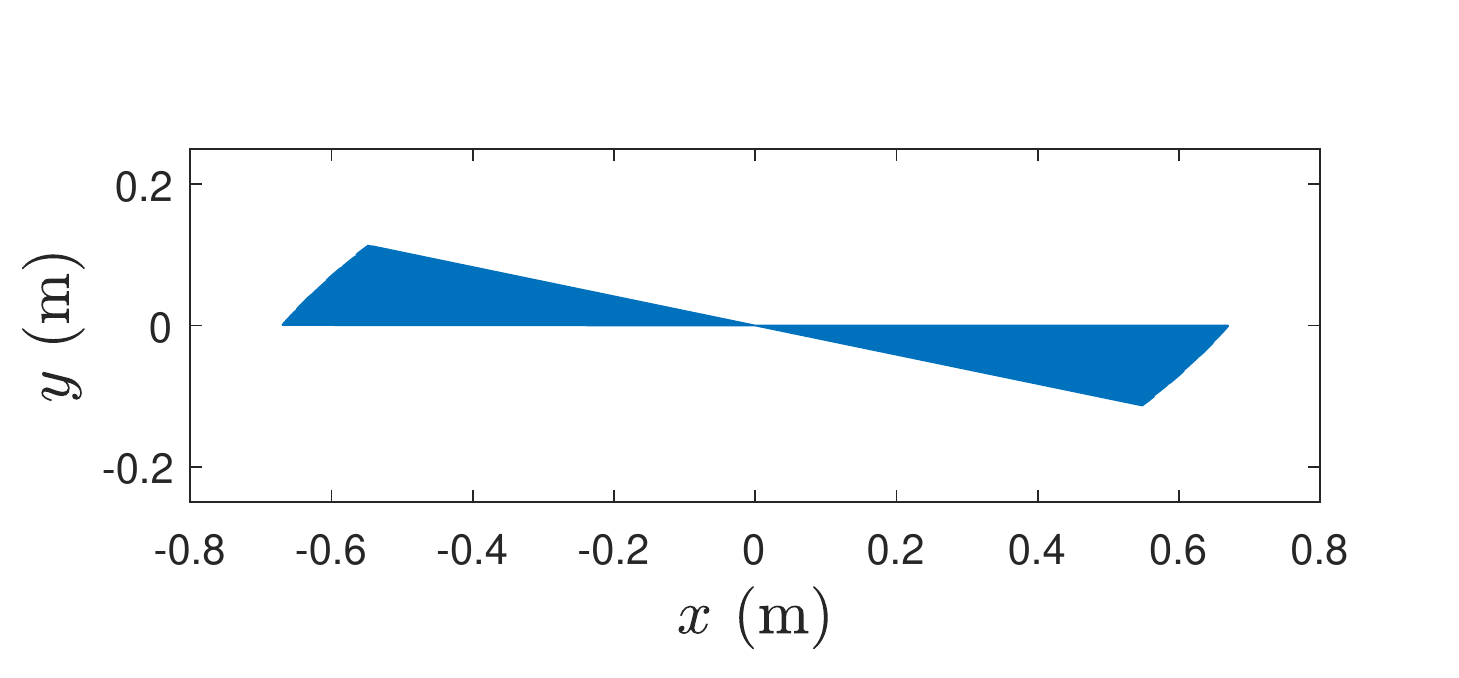}
    \caption{Contact (1st order)}
  \end{subfigure} \\
  \begin{subfigure}[b]{0.65\textwidth}
    \centering
    \includegraphics[width=\textwidth]{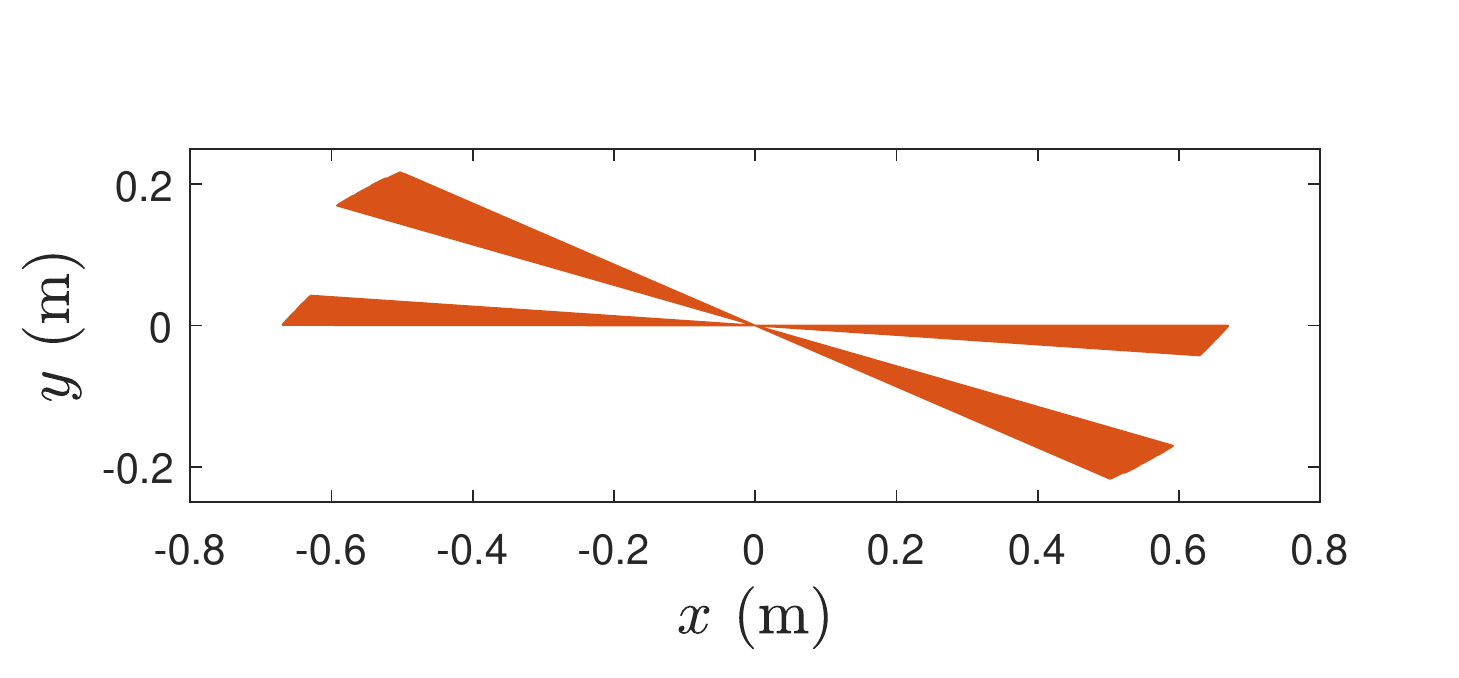}
    \caption{LA (1st order)}
  \end{subfigure}
  \begin{subfigure}[b]{0.65\textwidth}
    \centering
    \includegraphics[width=\textwidth]{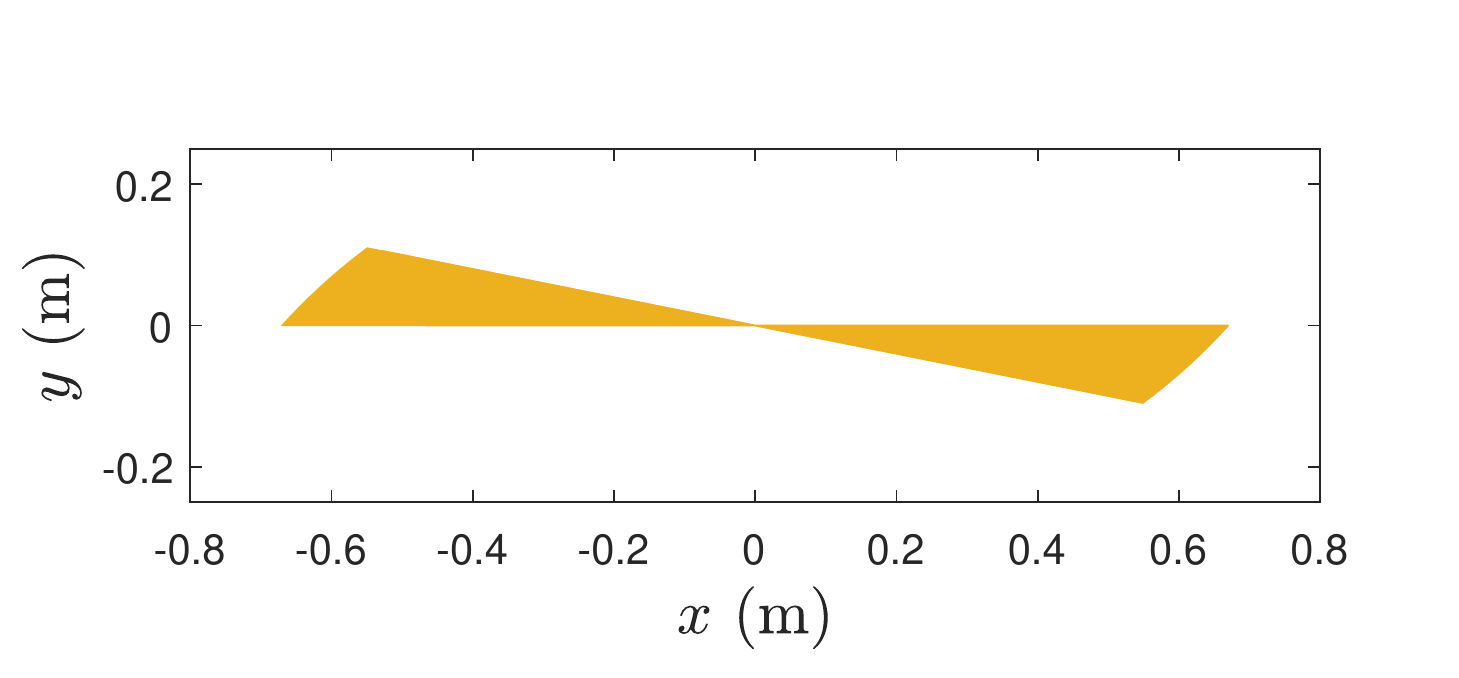}
    \caption{Reference (RKF45)}
  \end{subfigure}
  \caption{Foucault pendulum with initial conditions  $q(0) = (0,l/100)$, $\dot{q}(0) = (0,0)$. Dissipation parameter $\alpha = 1 \times 10^{-4}$. The graphics show the trajectories followed by the pendulum on the horizontal plane, for an oscillation time of 3600 s.}
  \label{fig:foucault_sim_2_oscillation_plane}
\end{figure}

\vskip 2mm
\noindent{\bf Experiment 2:} In order to simulate a more realistic situation, where dissipation is due to friction with the air, here we consider $\alpha = 1\times 10^{-4}$. In this setting, the contact integrator outcome is closer to the reference, both in terms of the error $\|q_{t_j}-q_{ref}(t_j)\|$ as well as in term of energy dissipation, as can be seen in Figure \ref{fig:foucault_sim_2_err_energy}. More interestingly, in this setting, the LA integrator displays anomalous behavior. Concretely the plane of oscillation varies discontinuously at a given time (see part b) of Figure \ref{fig:foucault_sim_2_oscillation_plane}. On the other hand, the contact integrator does not show this anomaly, as can be seen in part a) of Figure \ref{fig:foucault_sim_2_oscillation_plane}.

\subsection{Numerical Experiments: Falling rolling disk} For the falling rolling disk, we perform $4$ numerical experiments, in each case, using the MATLAB solver ode15i as the reference, with a step-size 10 times smaller than the one used for the proposed integrator. In all the four experiments, we use the following configuration for the system: $m = 5$ kg, $R = 0.5$ m, $I_A = \frac{1}{2}mR^2$ kgm$^2$, $I_T = \frac{1}{4}mR^2$ kgm$^2$, $g = 9.81$ m/s$^2$. The numerical parameters were set as: step-size $h = 0.1$, and tolerance $\epsilon = 1\times 10^{-6}$, both for the ode15i and for the modified multivariate Newton-Raphson method. The four sets of numerical experiments thus were used to explore different initial conditions and different values of the dissipation parameter $\alpha$. In each case, we use the MATLAB function ``decic'' to set consistent initial conditions for the solver ode15i.

\vskip 2mm
\noindent{\bf Experiment 1:} Here we consider a disk starting in a vertical position and an external force $F^\psi = \frac{1}{2}$ N which forces it to roll. Under these conditions, we perform two simulations, with dissipation parameter $\alpha = 0.005$ and $\alpha = 0.1$. The outcome of the proposed integrator and the reference are indistinguishable in these settings.

\vskip 2mm
\noindent{\bf Experiment 2:} For the second set of simulations we consider a disk starting in an inclined position, precisely, with $\theta_0 = \frac{\pi}{36}$ rad, and initial rolling velocity $\dot{\psi}_0 = 2\pi$. Under these conditions, we perform three simulations corresponding to $\alpha = 0$, $\alpha = 0.005$ and $\alpha = 0.1$. For the first two cases, the reference and the proposed integrator have indistinguishable behaviors. For the third case, we observe something interesting. In this case, the solver ode15i was not able to continue the simulation after $t\approx 11.1$ s, concretely displaying the following message. 

\noindent{\small \texttt{Warning: Failure at t=1.112335e+01.  Unable to meet integration tolerances\\
without reducing the step size below the smallest value allowed (3.951806e-14)\\
at time t}}.

\begin{figure}
  \centering
  \begin{subfigure}[b]{0.75\textwidth}
    \centering
    \includegraphics[width=\textwidth]{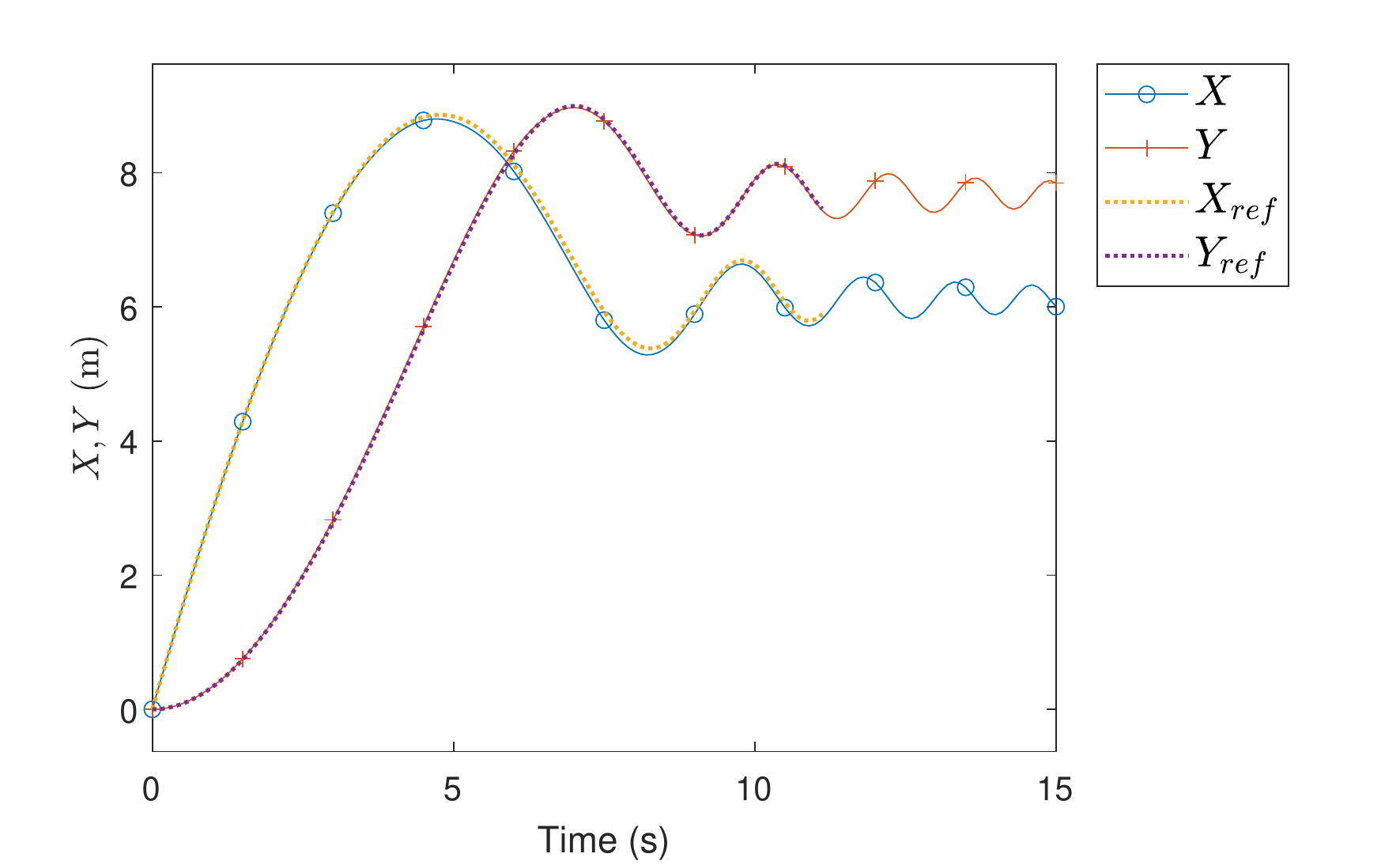}
  \end{subfigure}
  \par\bigskip
  \begin{subfigure}[b]{0.75\textwidth}
    \centering
    \includegraphics[width=\textwidth]{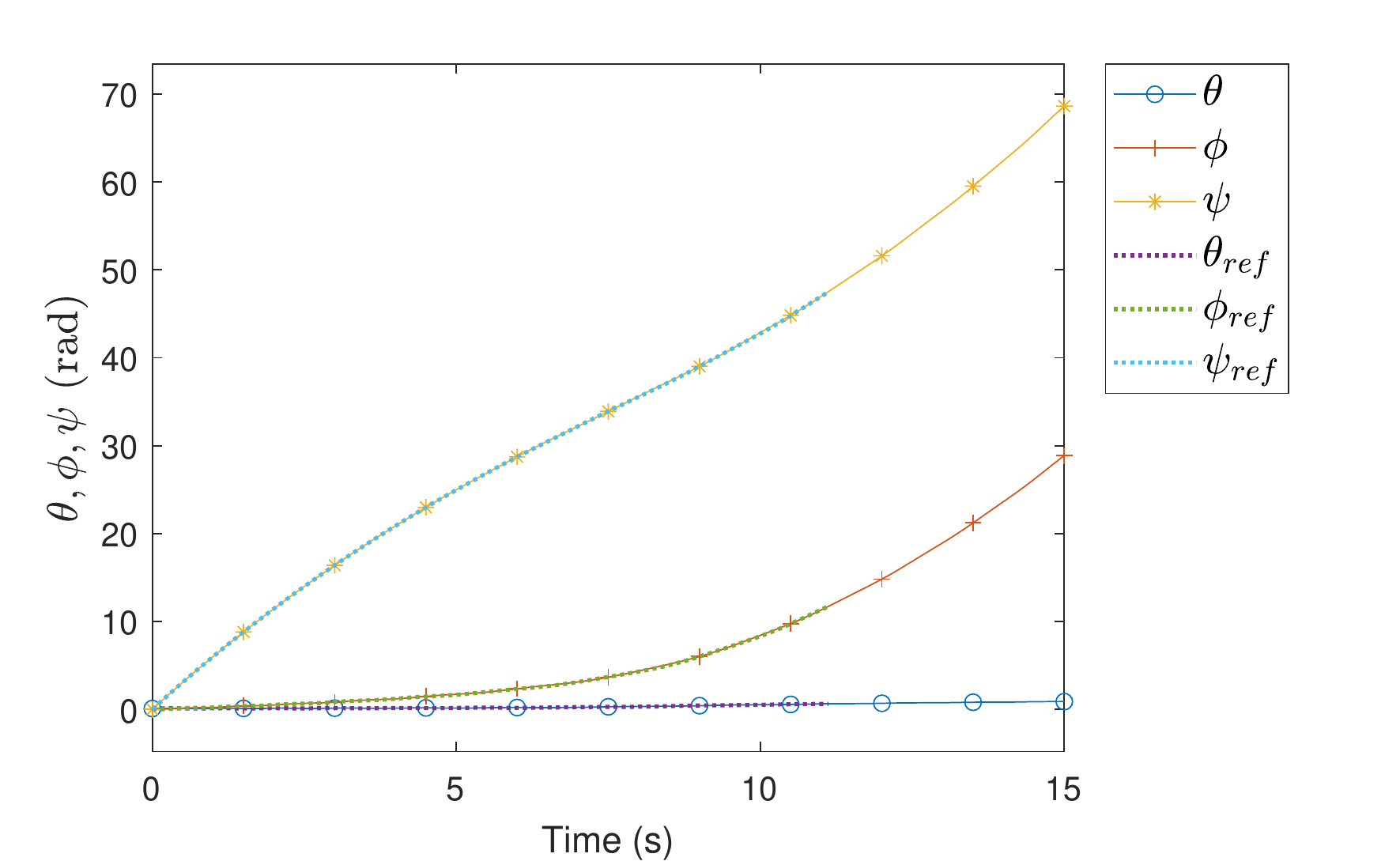}
  \end{subfigure}
  \caption{Experiment 2.3 of the falling disk with initial conditions ${q(0) = (0,0,\pi/36,0,0)}$, ${\dot{q}(0) = (\pi,0,0,0,2\pi)}$. Dissipation ${\alpha = 0.1}$ and ${F(t) = (0,0,0,0,0)}$. The graphics display the coordinate functions $(X,Y,\theta,\phi,\psi)$ of the system, obtained by the contact integrator and by ode15i, which is the reference method. For each coordinate, the $ref$ subindex indicates the solution given by the reference method.}
  \label{fig:disco_experimento_2_3_1}
\end{figure}

\begin{figure}
  \centering
  \begin{subfigure}[b]{0.75\textwidth}
    \centering
    \includegraphics[width=\textwidth]{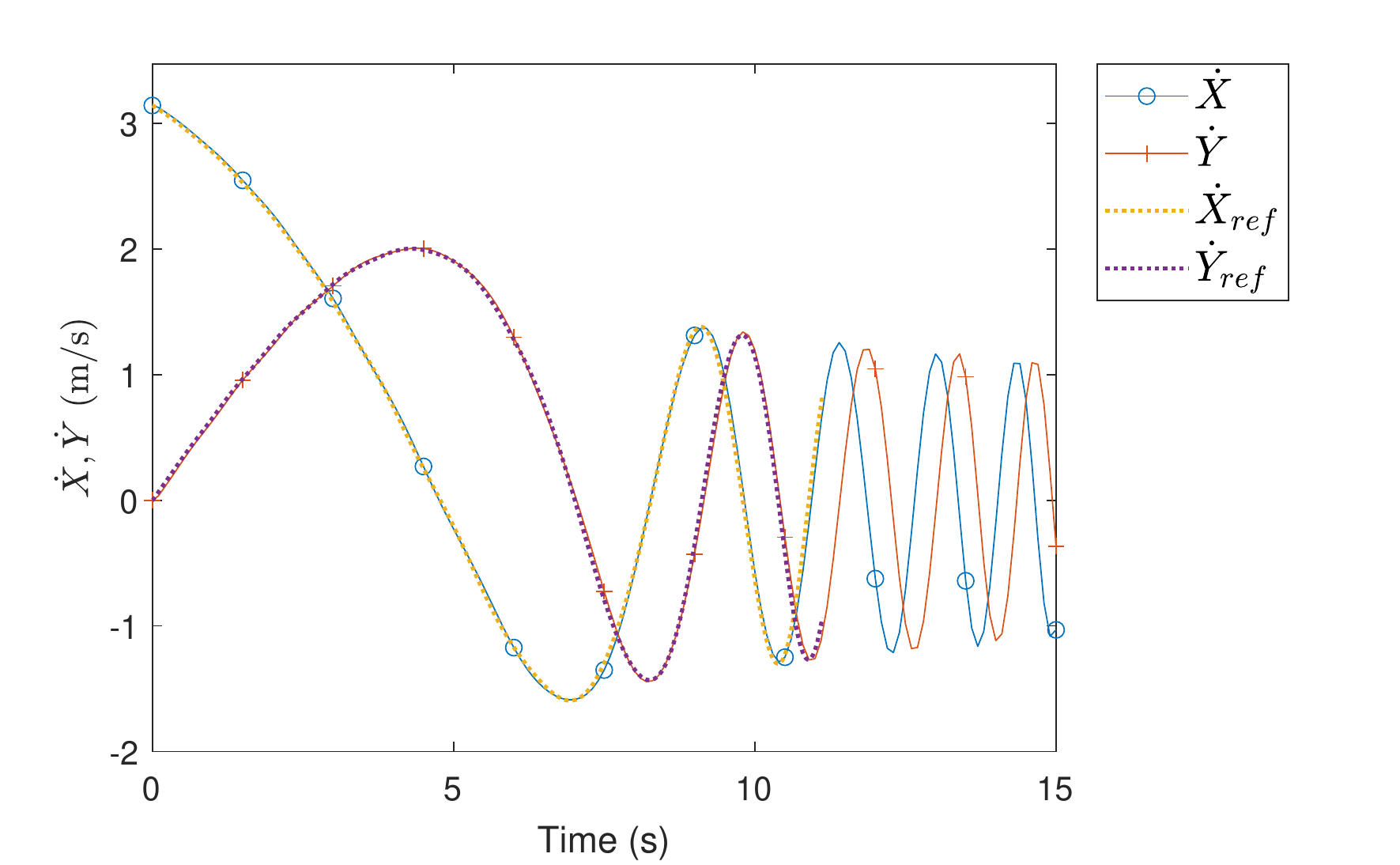}
  \end{subfigure}
  \par\bigskip
  \begin{subfigure}[b]{0.75\textwidth}
    \centering
    \includegraphics[width=\textwidth]{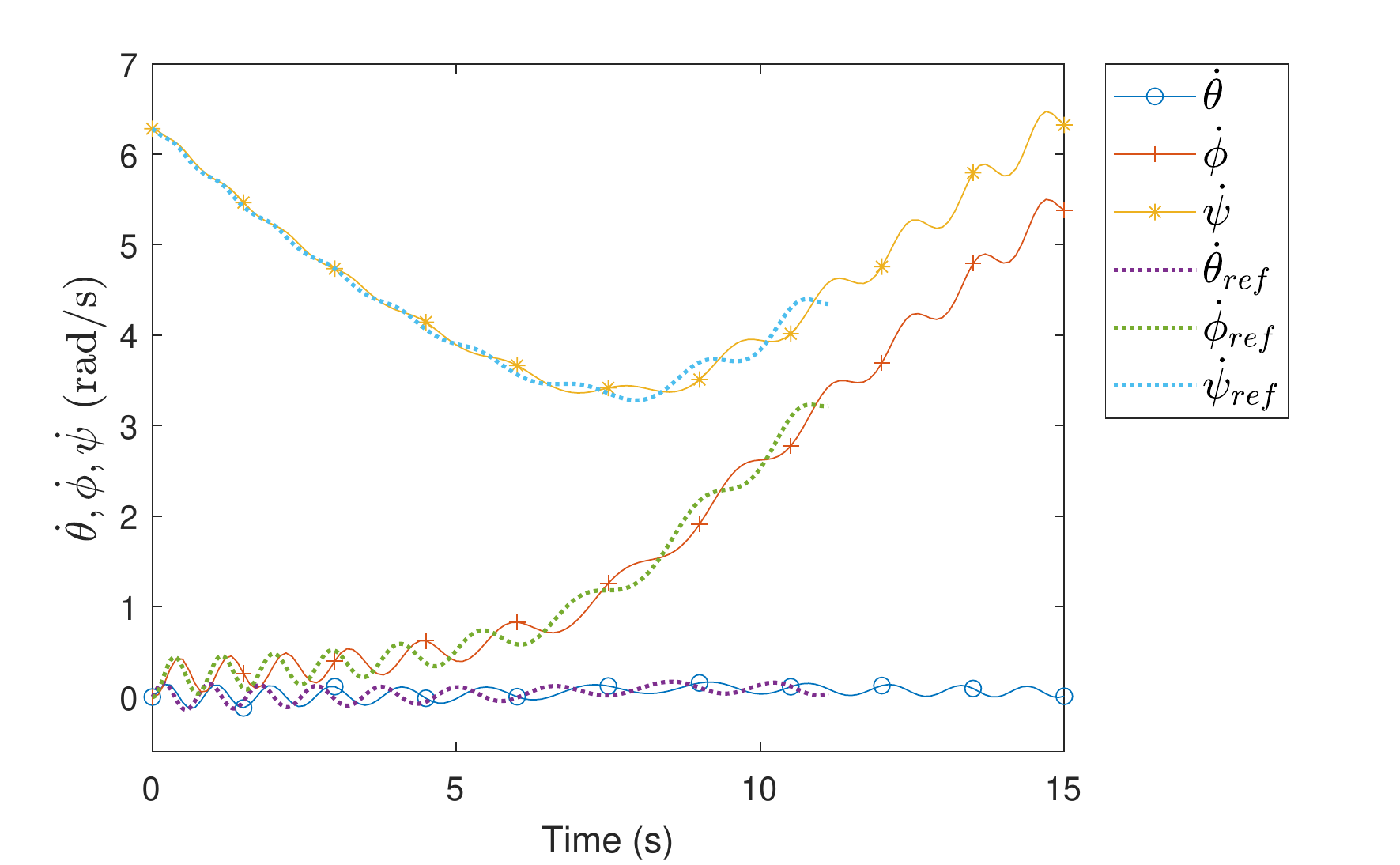}
  \end{subfigure}
  \caption{Experiment 2.3 of the falling disk with initial conditions ${q(0) = (0,0,\pi/36,0,0)}$, ${\dot{q}(0) = (\pi,0,0,0,2\pi)}$. Dissipation ${\alpha = 0.1}$ and ${F(t) = (0,0,0,0,0)}$. The graphics display the velocity functions $(\dot{X},\dot{Y},\dot{\theta},\dot{\phi},\dot{\psi})$ of the system, obtained by the contact integrator and by ode15i which is the reference method. For each velocity, the $ref$ subindex indicates the solution given by the reference method.}
  \label{fig:disco_experimento_2_3_2}
\end{figure}

\begin{figure}
  \centering
  \includegraphics[width=0.8\textwidth]{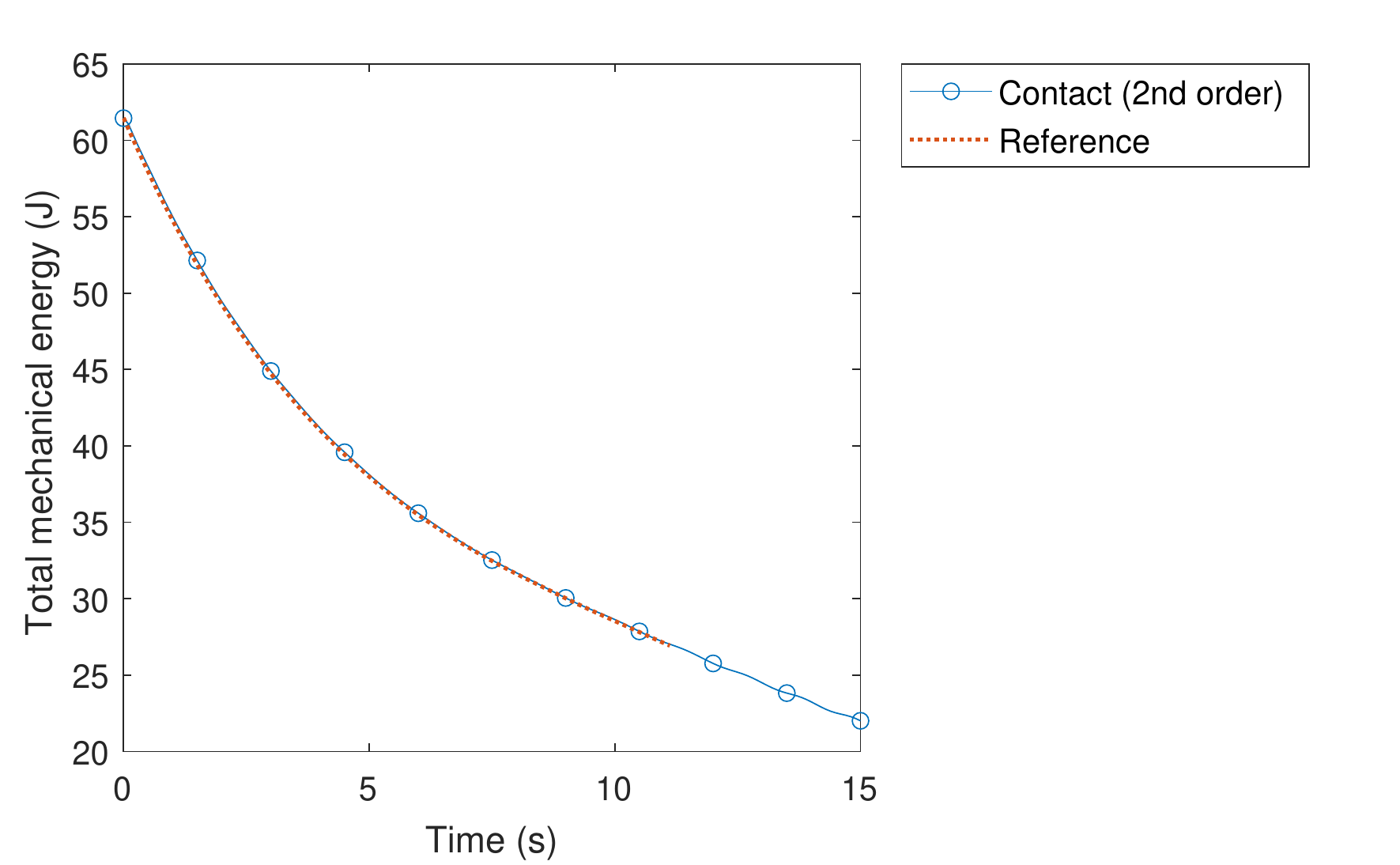}
  \caption{Experiment 2.3 of the falling disk with initial conditions ${q(0) = (0,0,\pi/36,0,0)}$, ${\dot{q}(0) = (\pi,0,0,0,2\pi)}$. Dissipation ${\alpha = 0.1}$ and ${F(t) = (0,0,0,0,0)}$. The graphic displays the energy functions obtained by the contact integrator and by ode15i which is the reference method.}
  \label{fig:disco_experimento_2_3_3}
\end{figure}

\begin{figure}
  \centering
  \includegraphics[width=0.8\textwidth]{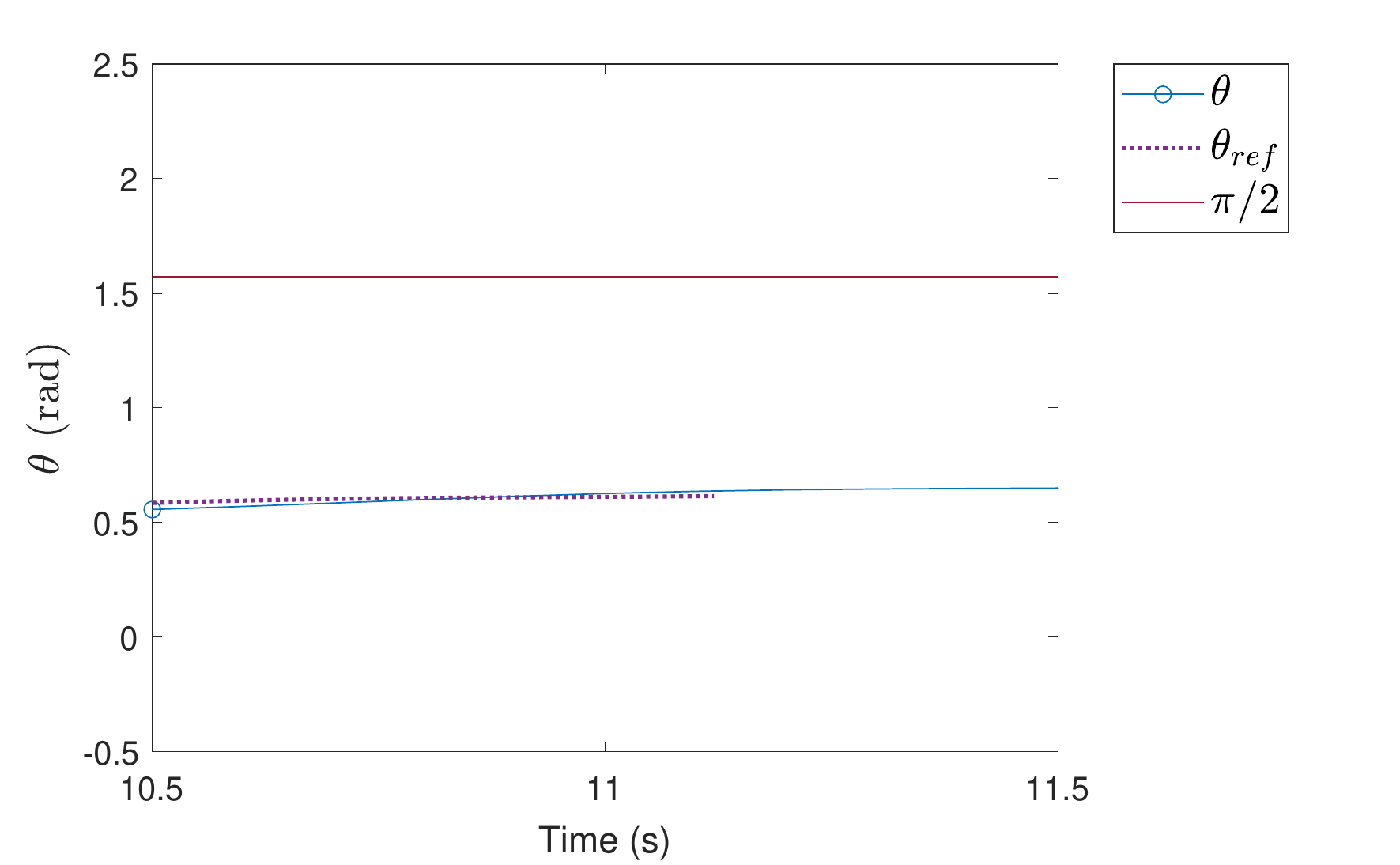}
  \caption{Experiment 2.3 of the falling disk with initial conditions ${q(0) = (0,0,\pi/36,0,0)}$, ${\dot{q}(0) = (\pi,0,0,0,2\pi)}$. Dissipation ${\alpha = 0.1}$ and ${F(t) = (0,0,0,0,0)}$. The graphic displays the function of the configuration variable $\theta$ obtained by the contact integrator and ode15i, which is the reference method. It can be observed that the disk does not fall completely for either of the two integrators in $t = 11.1$ s, as the disk falls completely when $\theta = \pi/2$.}
  \label{fig:disco_experimento_2_3_4}
\end{figure}

On the other hand, the proposed integrator was able to continue the simulation during the entire preset time interval. This phenomena can be observed in Figure \ref{fig:disco_experimento_2_3_1}, Figure \ref{fig:disco_experimento_2_3_2}, and Figure \ref{fig:disco_experimento_2_3_3}, displaying the evolution of the five generalized coordinates and the total energy of the system. It is worth noticing that the warning message displayed by the solver ode15i is not because the disk has reached a configuration corresponding to a completely fallen disk since, as can be seen in Figure \ref{fig:disco_experimento_2_3_4}, for the value of $t\approx 11.1$ s, the value of $\theta_{ref}$ is far from $\pi/2$, which is the value corresponding to a completely fallen disk. 

\begin{figure}
  \centering
  \begin{subfigure}[b]{0.75\textwidth}
    \centering
    \includegraphics[width=\textwidth]{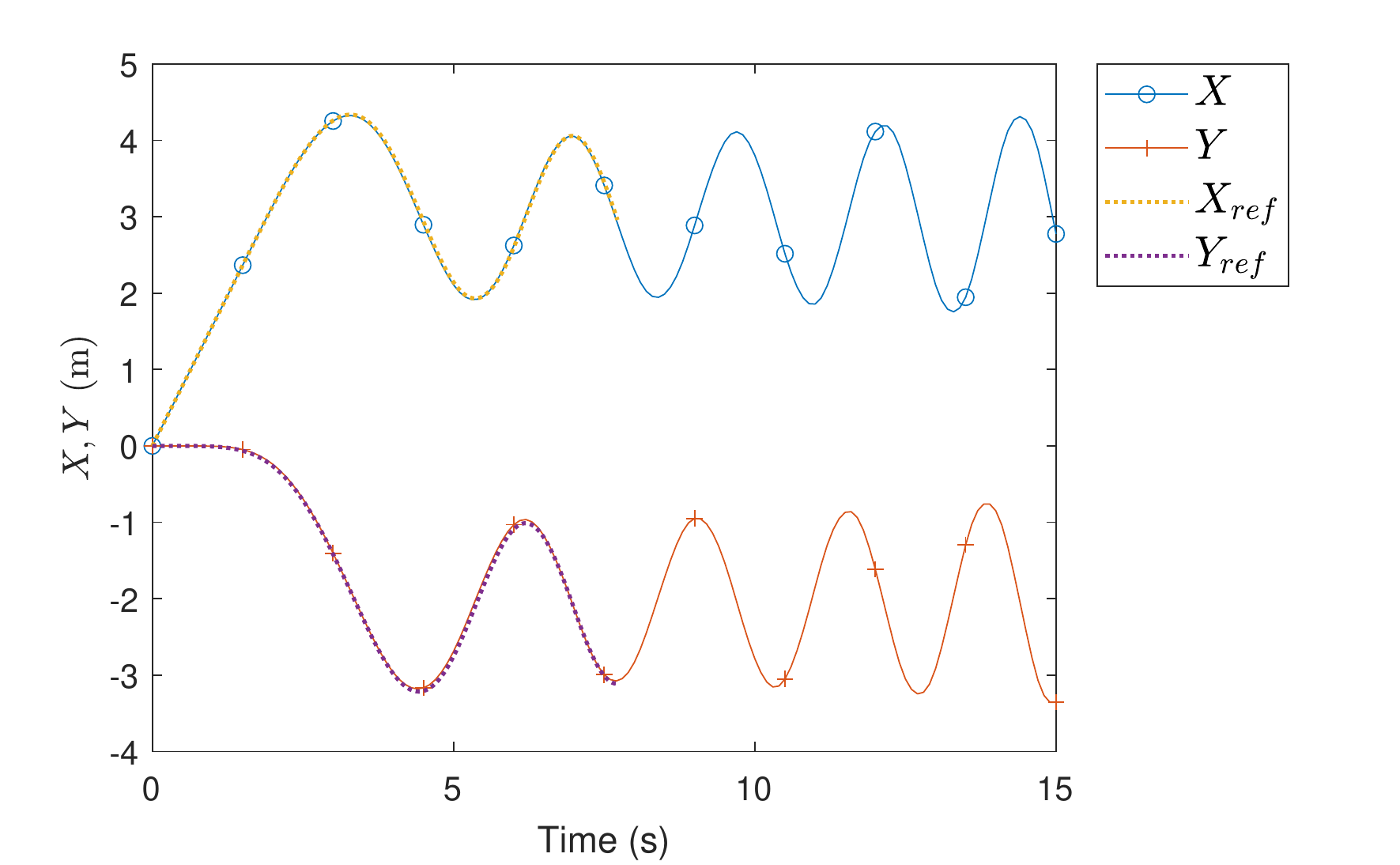}
  \end{subfigure}
  \par\bigskip
  \begin{subfigure}[b]{0.75\textwidth}
    \centering
    \includegraphics[width=\textwidth]{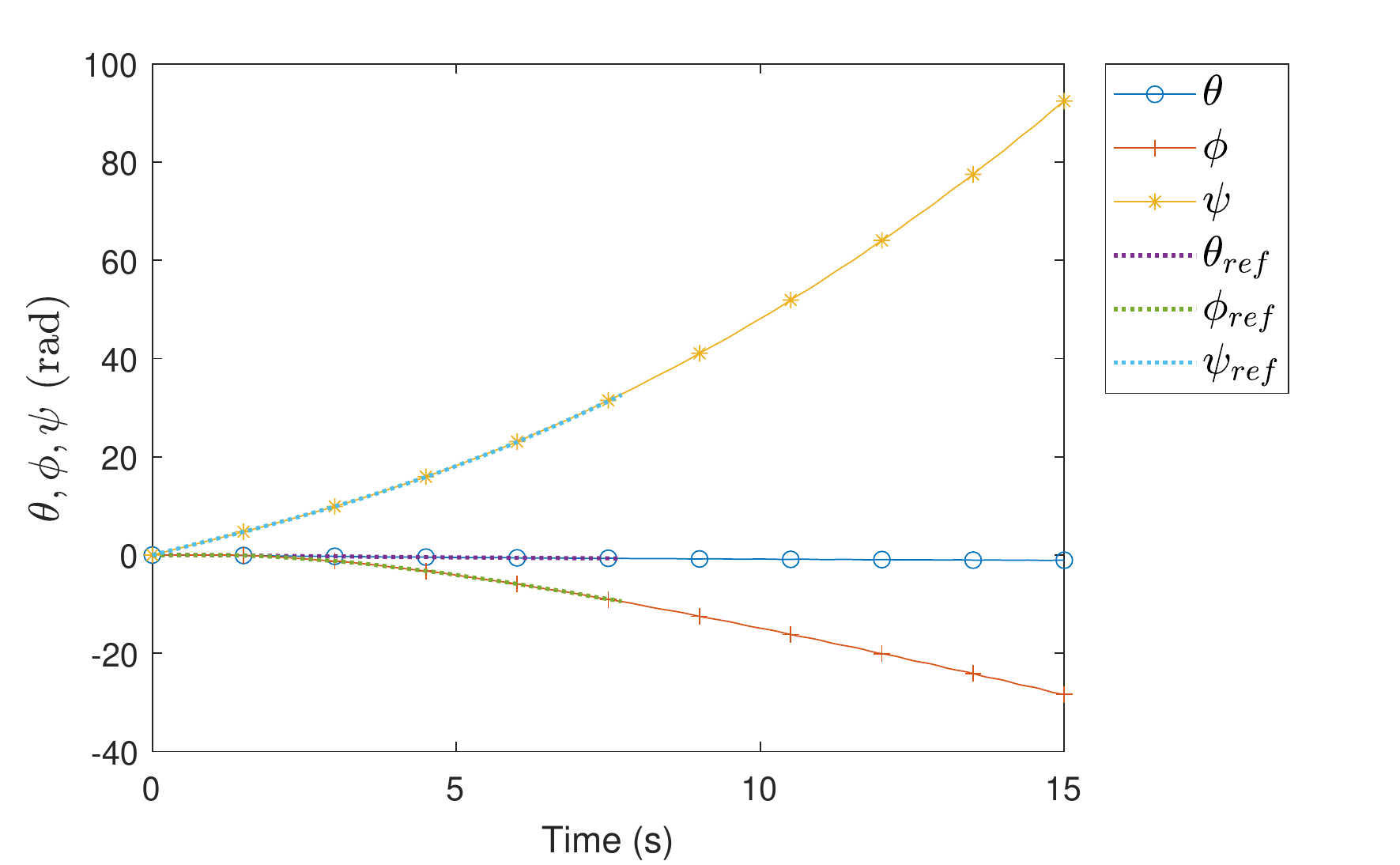}
  \end{subfigure}
  \caption{Experiment 3.1 of the falling rolling disk with initial condition $q(0)~=~(0,0,0,0,0)$, $\dot{q}(0) = (\pi/2,0,0,0,\pi)$. Dissipation $\alpha = 0$ and $F(t) = (0,0,0,t/16,t/16)$. The graphics show the time evolution of the generalized coordinates of the system for both integrators.}
  \label{fig:disco_experimento_3_1_1}
\end{figure}

\begin{figure}
  \centering
  \begin{subfigure}[b]{0.75\textwidth}
    \centering
    \includegraphics[width=\textwidth]{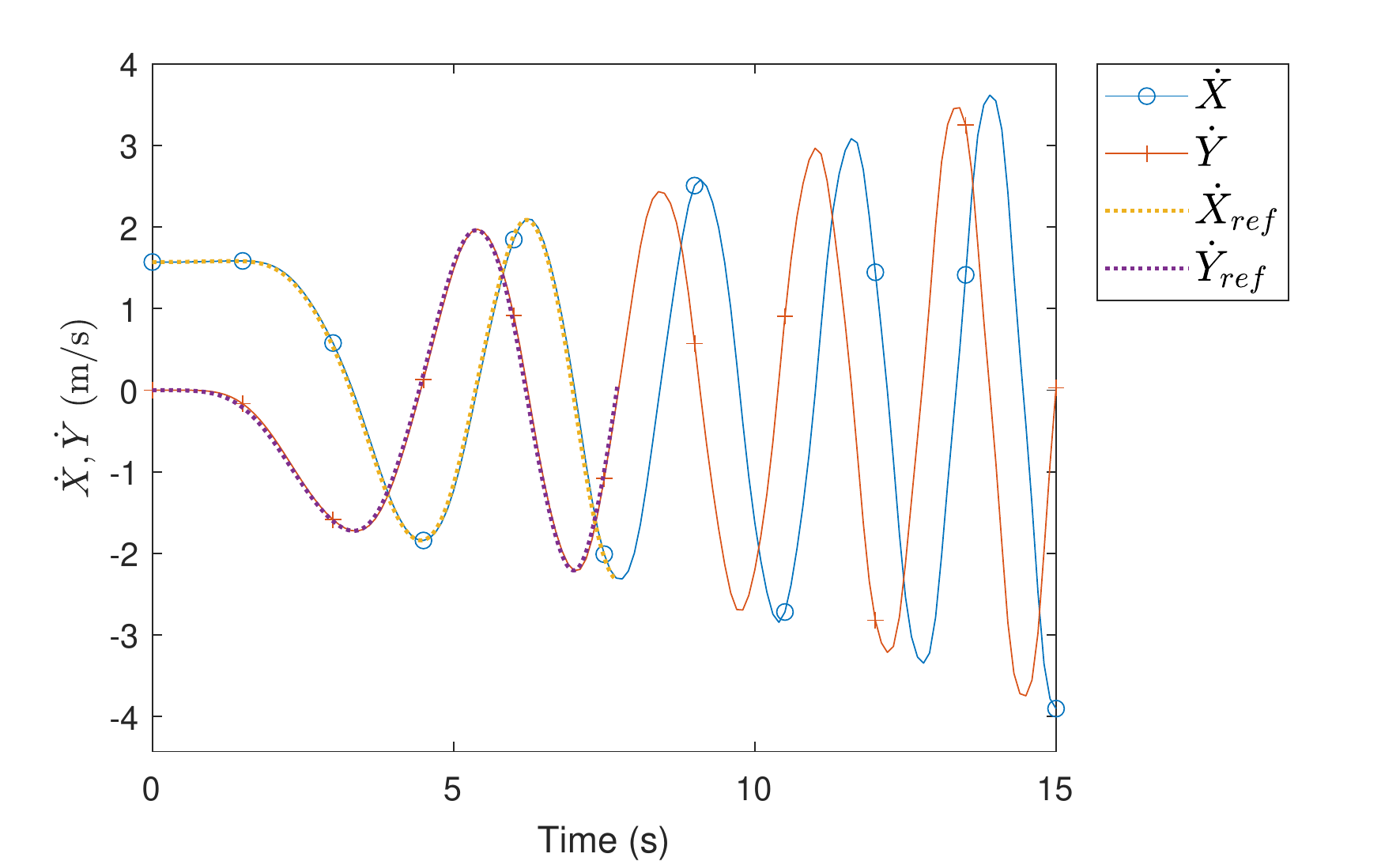}
  \end{subfigure}
  \par\bigskip
  \begin{subfigure}[b]{0.75\textwidth}
    \centering
    \includegraphics[width=\textwidth]{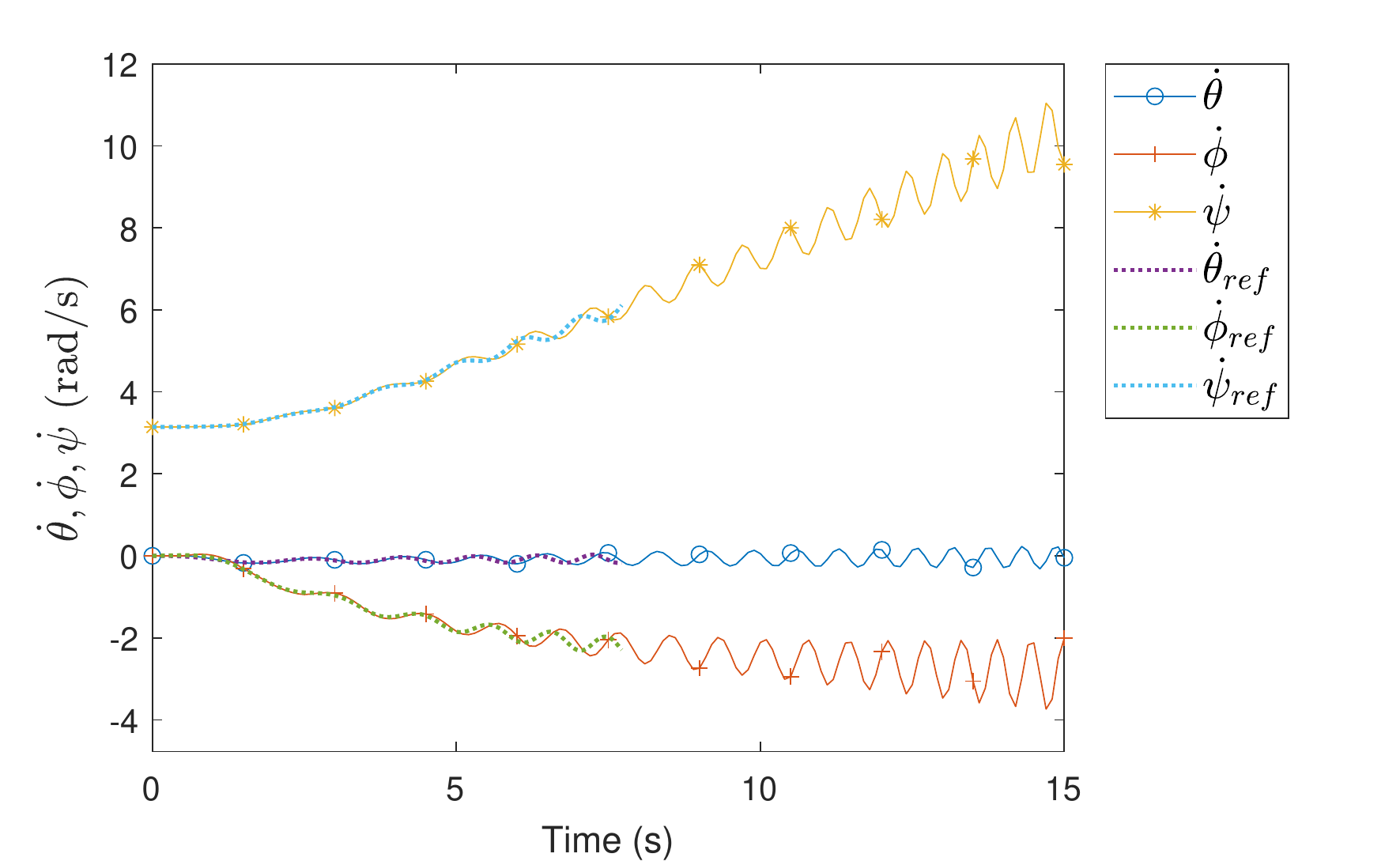}
  \end{subfigure}
  \caption{Experiment 3.1 of the falling rolling disk with initial condition $q(0)~=~(0,0,0,0,0)$, $\dot{q}(0) = (\pi/2,0,0,0,\pi)$. Dissipation $\alpha = 0$ and $F(t) = (0,0,0,t/16,t/16)$. The graphics show the time evolution of the generalized velocities of the system for both integrators.}
  \label{fig:disco_experimento_3_1_2}
\end{figure}

\begin{figure}
  \centering
  \includegraphics[width=0.8\textwidth]{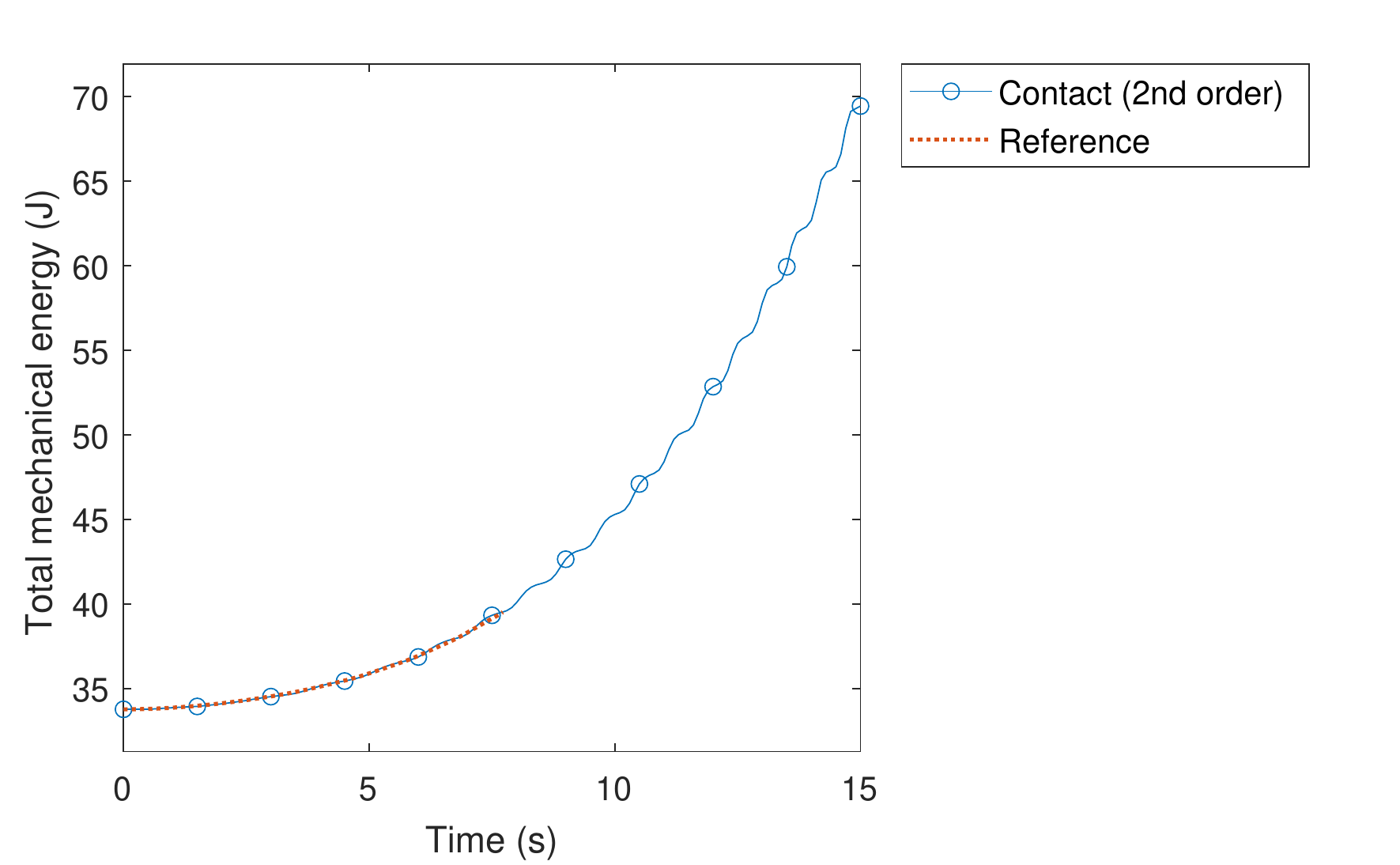}
  \caption{Experiment 3.1 of the falling rolling disk with initial condition $q(0)~=~(0,0,0,0,0)$, $\dot{q}(0) = (\pi/2,0,0,0,\pi)$. Dissipation $\alpha = 0$ and $F(t) = (0,0,0,t/16,t/16)$. The graphic shows the energy function for both integrators.}
  \label{fig:disco_experimento_3_1_3}
\end{figure}

\begin{figure}
  \centering
  \includegraphics[width=0.8\textwidth]{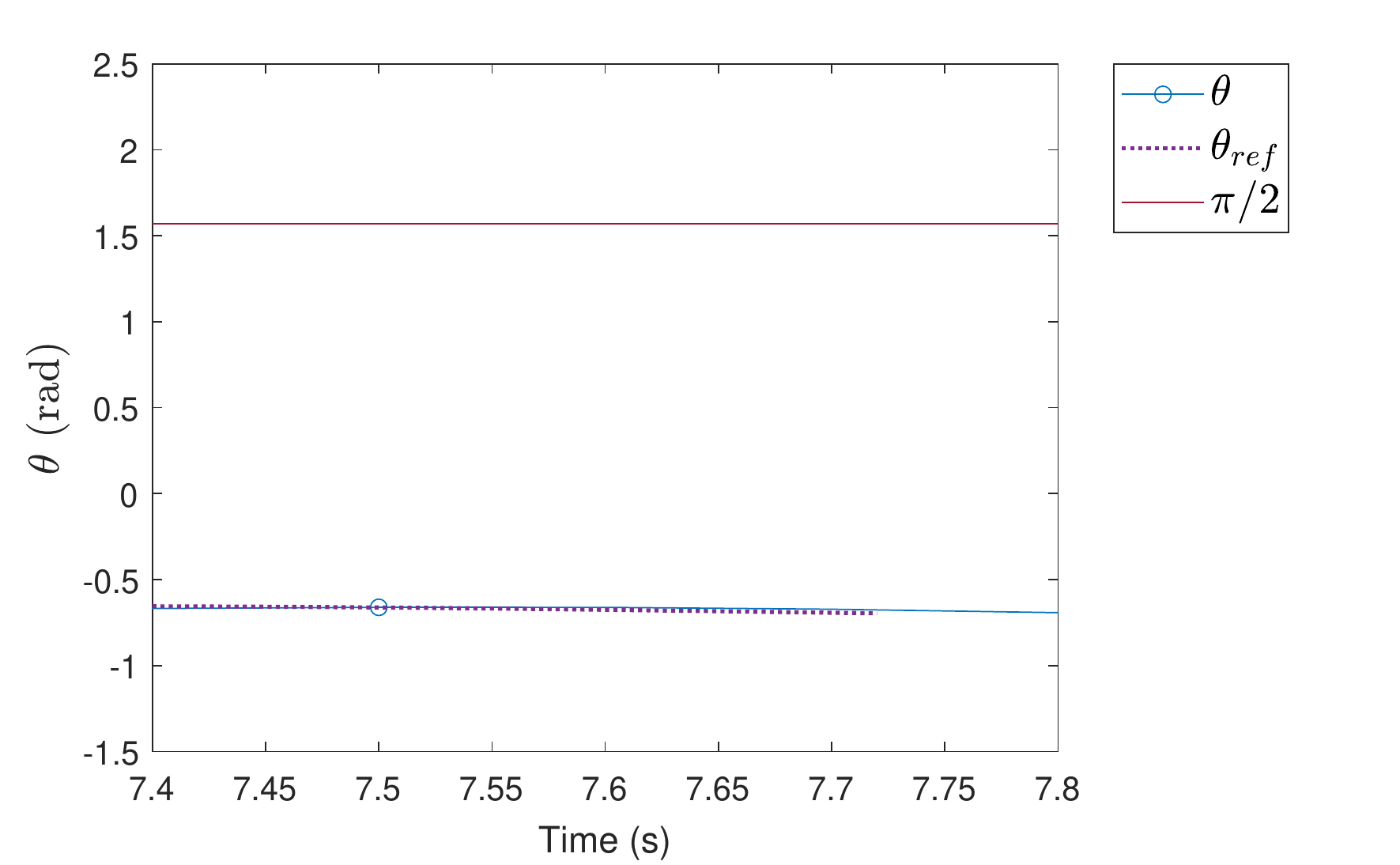}
  \caption{Experiment 3.1 of the falling rolling disk with initial condition $q(0)~=~(0,0,0,0,0)$, $\dot{q}(0) = (\pi/2,0,0,0,\pi)$. Dissipation $\alpha = 0$ and $F(t) = (0,0,0,t/16,t/16)$. The graphic displays the function of the configuration variable $\theta$ obtained by the contact integrator and ode15i, which is the reference method. It can be observed that the disk does not fall completely for either of the two integrators in $t \approx 7.7$ s, as the disk falls completely when $\theta = \pi/2$.}
  \label{fig:disco_experimento_3_1_4}
\end{figure}

\vskip 2mm
\noindent{\bf Experiment 3:} For this set of simulations, we consider again a disk starting in a vertical position, with initial rolling velocity $\dot{\psi}_0 = \pi$. Additionally, here we consider two applied forces $F^\psi ; \frac{t}{16}$ N and $F^\phi = \frac{t}{16}$ N, changing the rolling speed and the direction, respectively. We perform simulations for three values of the dissipation parameter, namely: $\alpha = 0$, $\alpha = 0.005$ and $\alpha = 0.1$.
The interesting case now corresponds to $\alpha = 0$, where once again, the reference solver ode15i fails to complete the simulation, displaying the same warning message as before for $t \approx 7.7$ s. As before, this failure is not due to a configuration corresponding to a completely fallen disk, as can be seen in Figure \ref{fig:disco_experimento_3_1_4}. On the other hand, the proposed integrator does not suffer from this issue, as can be observed in Figure \ref{fig:disco_experimento_3_1_1}, Figure \ref{fig:disco_experimento_3_1_2}, and Figure \ref{fig:disco_experimento_3_1_3}. For the other values of $\alpha$, the outcomes of both integrators are qualitatively well-behaved.

\begin{figure}
  \centering
  \begin{subfigure}[b]{0.75\textwidth}
    \centering
    \includegraphics[width=\textwidth]{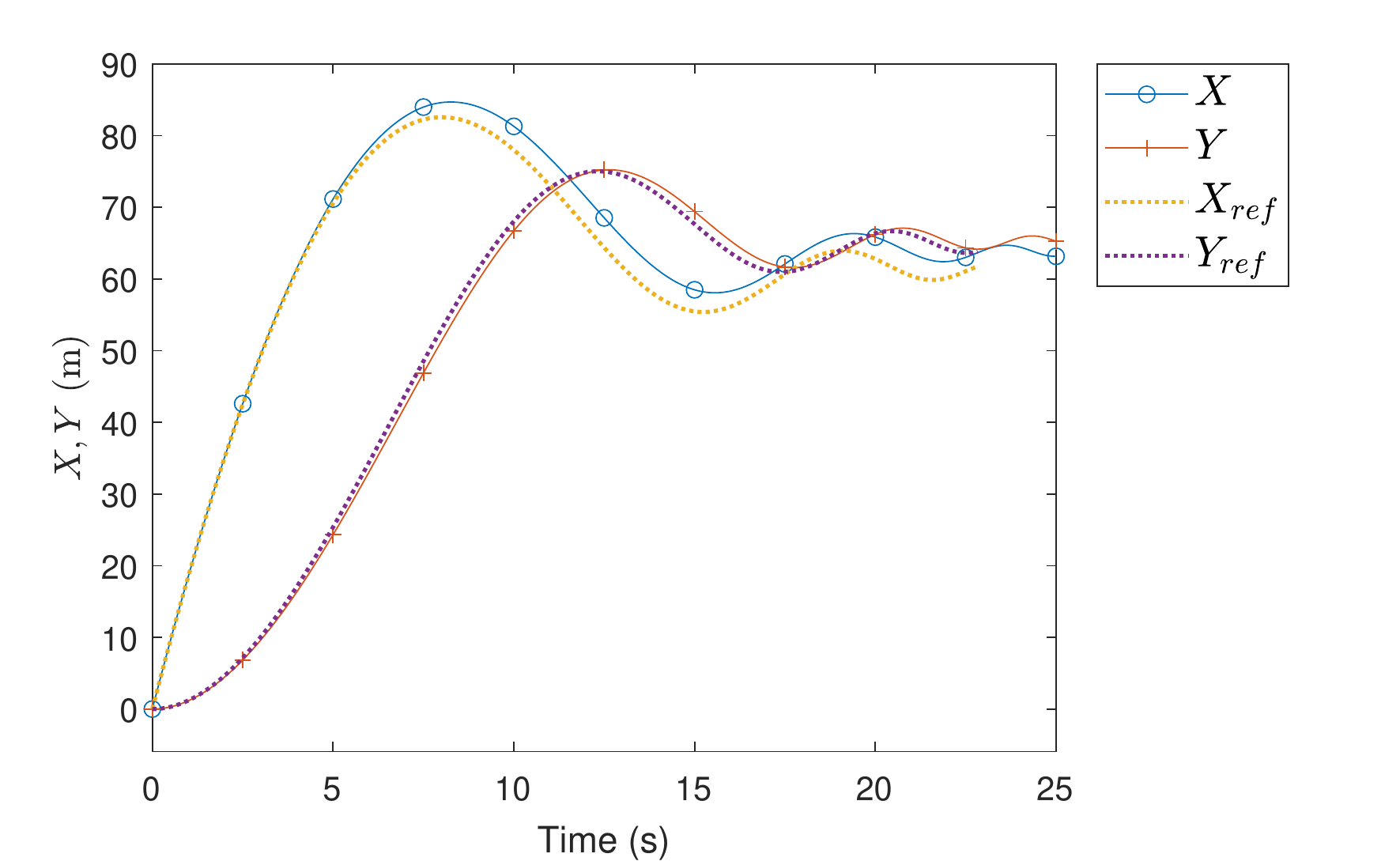}
  \end{subfigure}
  \par\bigskip
  \begin{subfigure}[b]{0.75\textwidth}
    \centering
    \includegraphics[width=\textwidth]{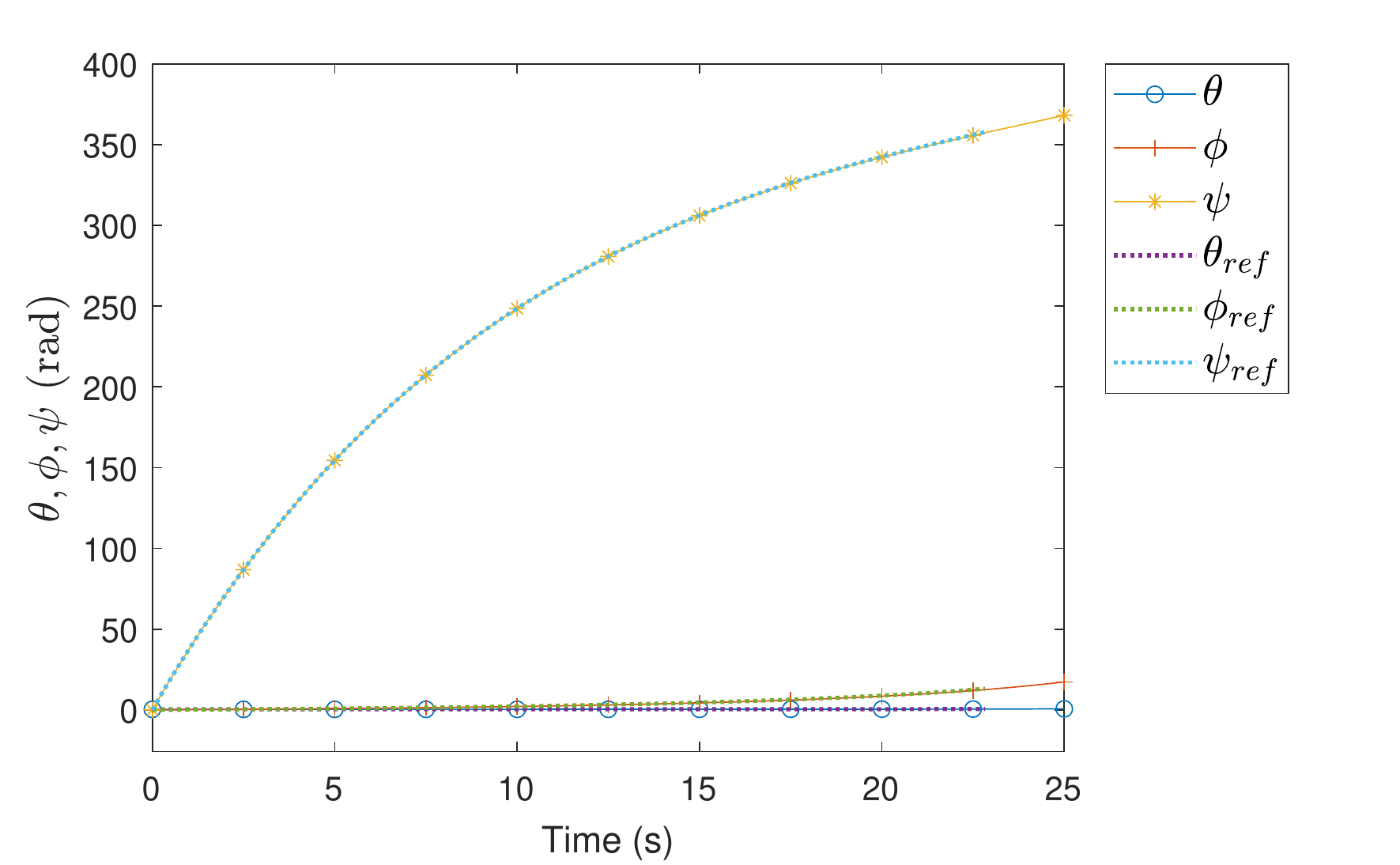}
  \end{subfigure}
  \caption{Experiment 4.3 of the falling rolling disk with initial conditions $q(0)~=~(0,0,20\pi/180,0,0)$, $\dot{q}(0) = (\pi/2,0,0,-3\pi/10,\dot{\psi}_0)$, with $\dot{\psi}_0 = ((I_T - I_A - m R^2) \sin(\theta_0) \dot{\phi}_0^2 - m g R)/((I_A + m R^2) \tan(\theta_0) \dot{\phi}_0)$. Dissipation $\alpha = 0.1$ and $F(t) = (0,0,0,0,0)$. The graphics show the time evolution of the generalized coordinates for both integrators.}
  \label{fig:disco_experimento_4_3_1}
\end{figure}

\begin{figure}
  \centering
  \begin{subfigure}[b]{0.75\textwidth}
    \centering
    \includegraphics[width=\textwidth]{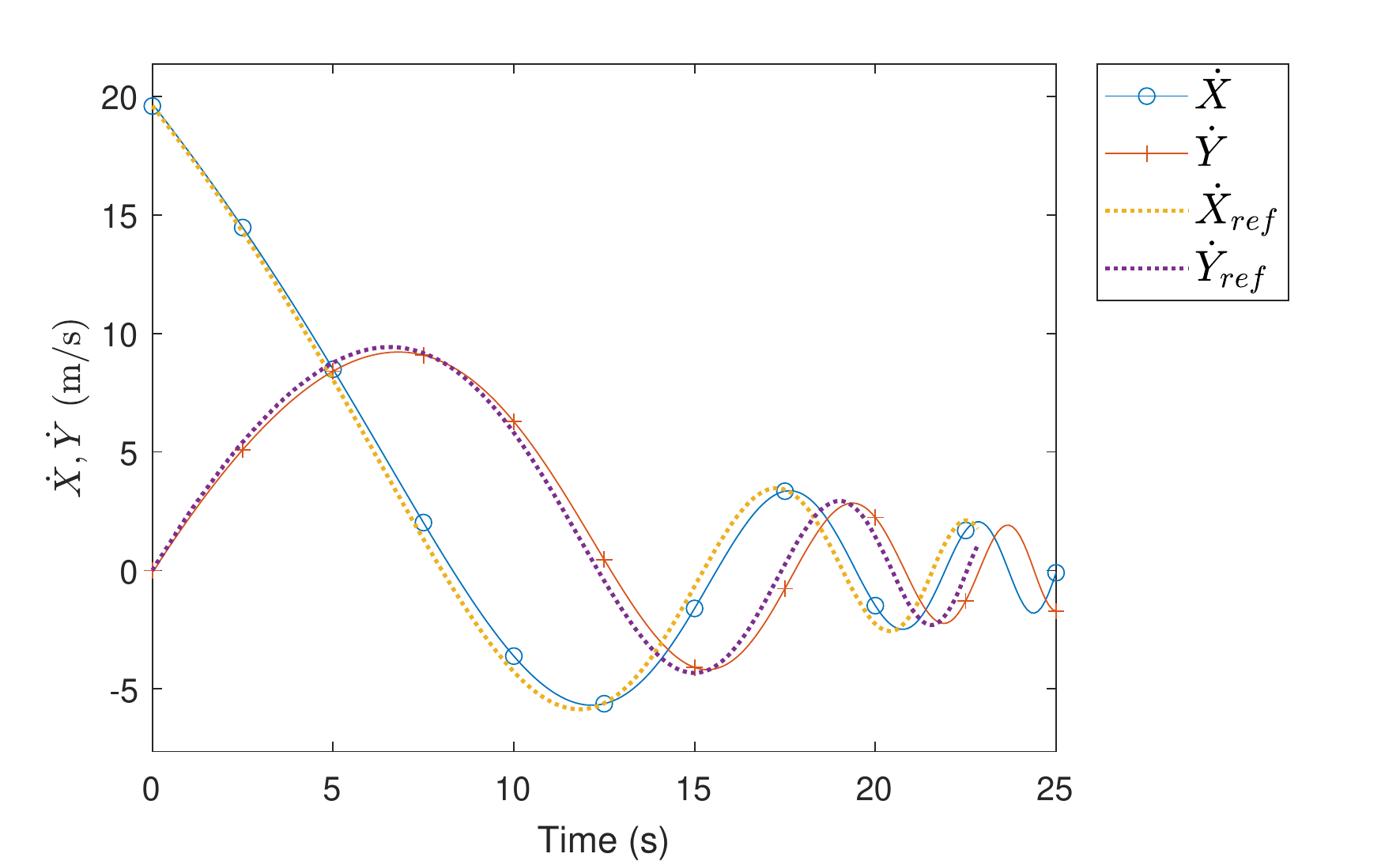}
  \end{subfigure}
  \par\bigskip
  \begin{subfigure}[b]{0.75\textwidth}
    \centering
    \includegraphics[width=\textwidth]{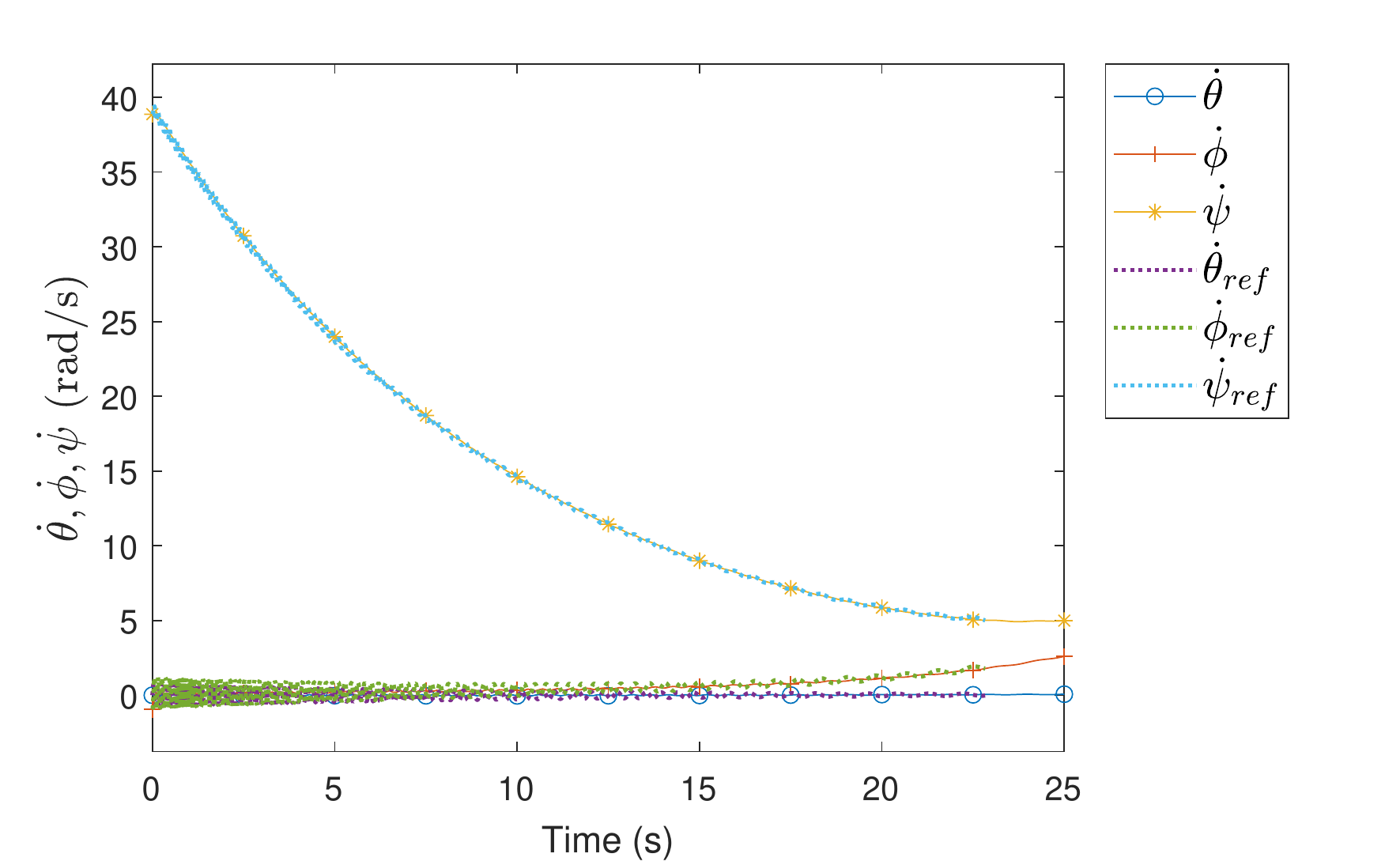}
  \end{subfigure}
  \caption{Experiment 4.3 of the falling rolling disk with initial conditions $q(0)~=~(0,0,20\pi/180,0,0)$, $\dot{q}(0) = (\pi/2,0,0,-3\pi/10,\dot{\psi}_0)$, with $\dot{\psi}_0 = ((I_T - I_A - m R^2) \sin(\theta_0) \dot{\phi}_0^2 - m g R)/((I_A + m R^2) \tan(\theta_0) \dot{\phi}_0)$. Dissipation $\alpha = 0.1$ and $F(t) = (0,0,0,0,0)$. The graphics show the time evolution of the generalized velocities for both integrators.}
  \label{fig:disco_experimento_4_3_2}
\end{figure}

\begin{figure}
  \centering
  \includegraphics[width=0.8\textwidth]{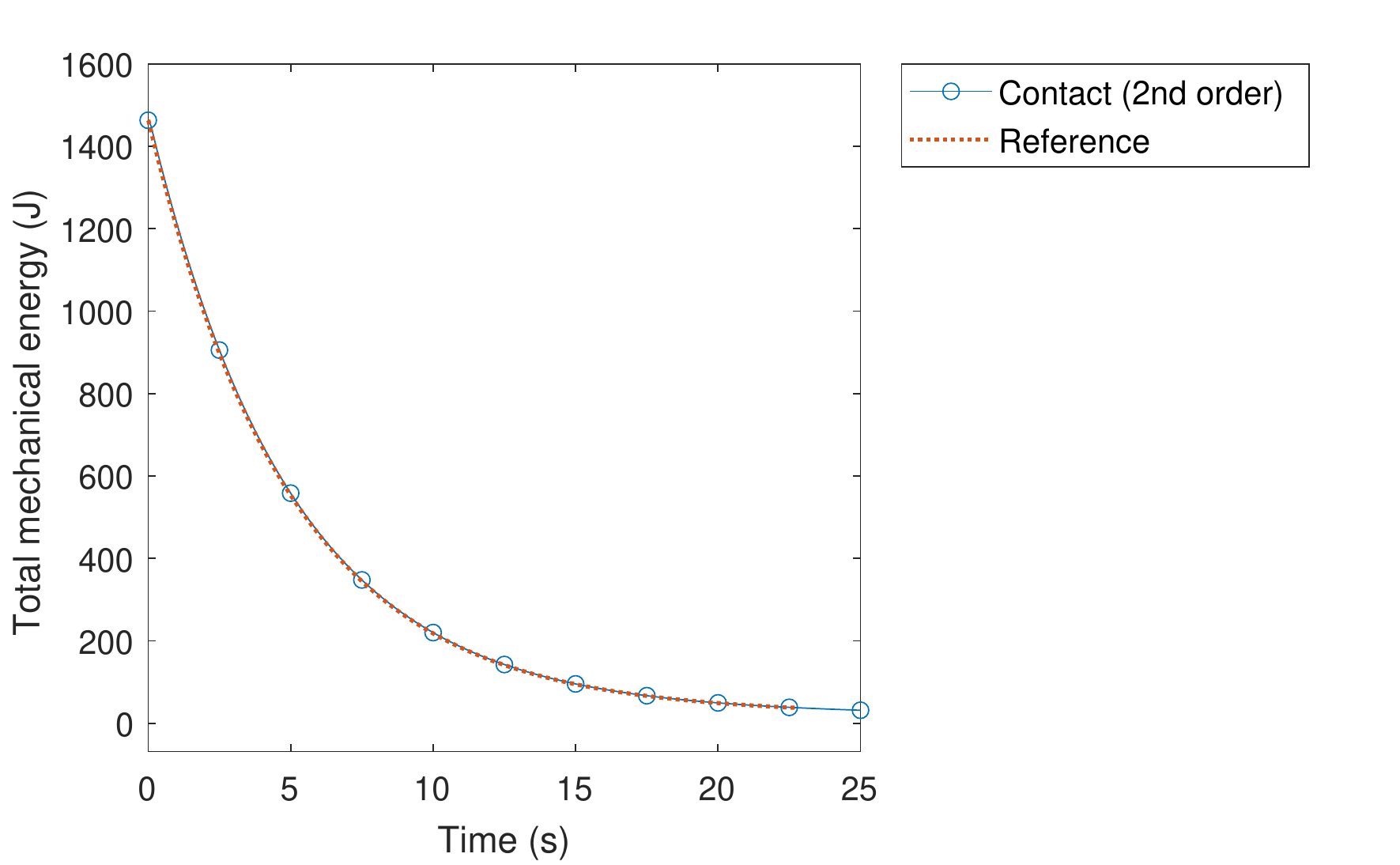}
  \caption{Experiment 4.3 of the falling rolling disk with initial conditions $q(0)~=~(0,0,20\pi/180,0,0)$, $\dot{q}(0) = (\pi/2,0,0,-3\pi/10,\dot{\psi}_0)$, with $\dot{\psi}_0 = ((I_T - I_A - m R^2) \sin(\theta_0) \dot{\phi}_0^2 - m g R)/((I_A + m R^2) \tan(\theta_0) \dot{\phi}_0)$. Dissipation $\alpha = 0.1$ and $F(t) = (0,0,0,0,0)$. The graphic shows the time evolution of the energy of the system for both integrators.}
  \label{fig:disco_experimento_4_3_3}
\end{figure}

\begin{figure}
  \centering
  \includegraphics[width=0.8\textwidth]{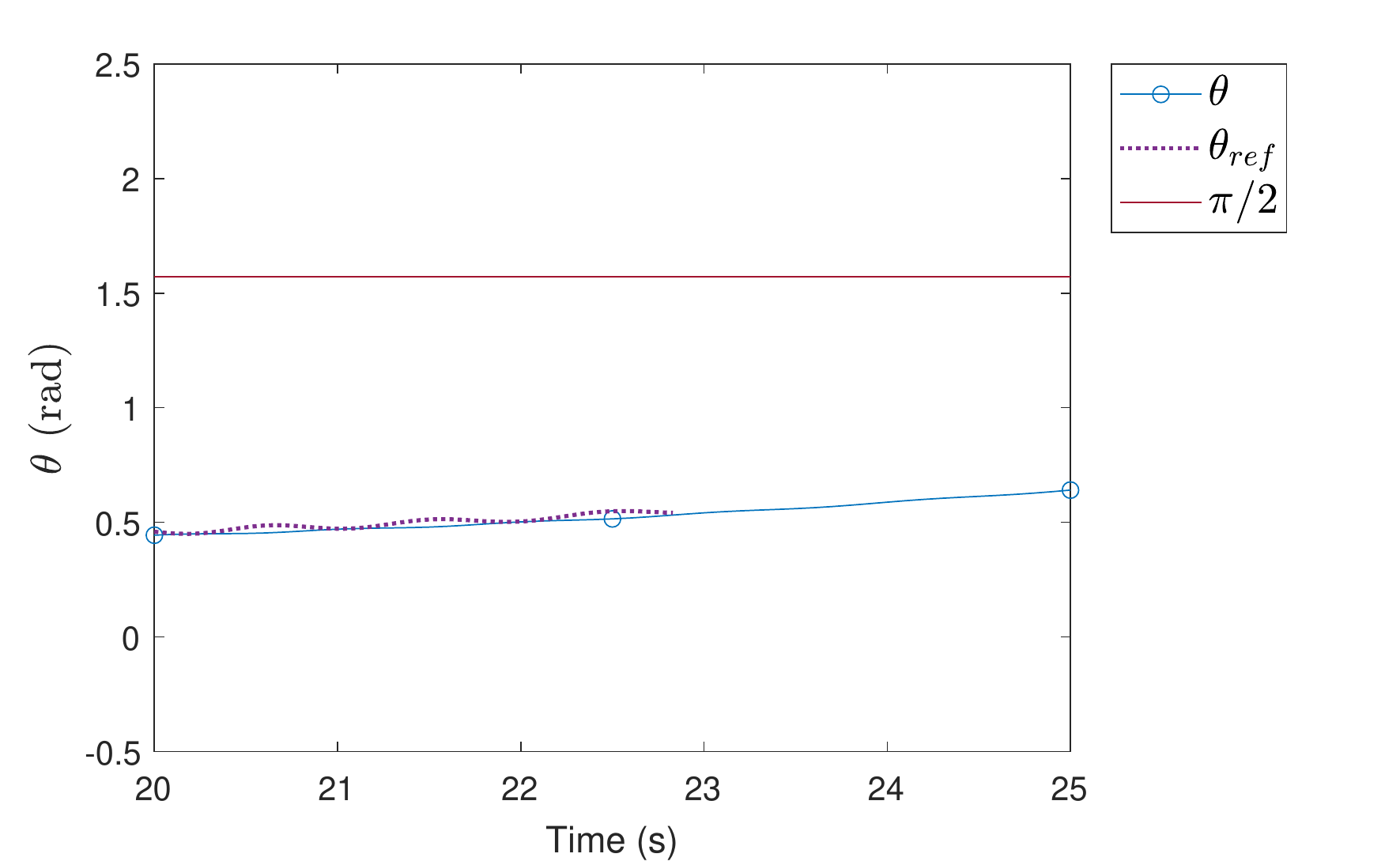}
  \caption{Experiment 4.3 of the falling rolling disk with initial conditions $q(0)~=~(0,0,20\pi/180,0,0)$, $\dot{q}(0) = (\pi/2,0,0,-3\pi/10,\dot{\psi}_0)$, with $\dot{\psi}_0 = ((I_T - I_A - m R^2) \sin(\theta_0) \dot{\phi}_0^2 - m g R)/((I_A + m R^2) \tan(\theta_0) \dot{\phi}_0)$. Dissipation $\alpha = 0.1$ and $F(t) = (0,0,0,0,0)$. The graphic displays the function of the configuration variable $\theta$ obtained by the contact integrator and ode15i, which is the reference method. It can be observed that the disk does not fall completely for either of the two integrators in $t \approx 22.8$ s, as the disk falls completely when $\theta = \pi/2$.}
  \label{fig:disco_experimento_4_3_4}
\end{figure}

\vskip 2mm
\noindent{\bf Experiment 4:} This set of simulations considers a condition under which the disk follows a stable circular path. This evolution is possible if the inclination angle $\theta = \theta_0$, the precession $\dot{\phi} = \dot{\phi}_0$, and the rolling velocity $\dot{\psi} = \dot{\psi}_0$ satisfy the following restriction
\[(I_T-I_A-mR^2)^2\sin(\theta_0)\dot{\theta}_0^2 - (I_A+mR^2)\tan(\theta_0)\dot{\theta}_0\dot{\psi}_0 - mgR = 0.\]
We perform again three simulations, corresponding to $\alpha = 0$, $\alpha =0.005$ and $\alpha = 0.1$. In these simulations, both integrators display, qualitatively,  the same outcome, but once again, for $\alpha = 0.1$, the solver ode15i was unable to complete the simulation, displaying the same warning message as in the previous situations, while the proposed integrator completed the simulation satisfactorily. We summarize this observations with Figures \ref{fig:disco_experimento_4_3_1} to \ref{fig:disco_experimento_4_3_4}.

\section{Concluding remarks}\label{sec:Conclusions}

In this work, we propose an integrator for nonholonomic, nonconservative mechanical systems. To do so, we combine previous works dealing with a derivation of a geometric integrator based on Herglotz variational principle, as well as one incorporating nonholonomic constraints for integrators derived from the Lagrange-d'Alembert principle.  {For many interesting mechanical systems, the use of Herglotz principle allows to incorporate dissipation phenomena directly into the contact-type Lagrangian, providing a more intrinsic treatment of these kinds of mechanical systems.}

We validate our integrator by numerical experiments performed over two mechanical systems, namely the Foucault pendulum with dissipation and the rolling falling disk with dissipation and a forcing term. The validation was done by comparing the numerical outcome of our integrator with standard methods such as the fourth-order Runge-Kutta-Fehlberg (for the Foucault pendulum) and the ode15i solver available in Matlab (for the rolling falling disk). Also, for the Foucault pendulum, we contrast our integrator against the one coming from the Lagrange-d'Alembert principle, which is somehow the ``standard'' variational-like approach to nonholonomic systems.

The numerical experiments show that our integrator has good qualitative behavior, matching the outcome of the references and, in some situations, outperforming them, at least in the qualitative sense. This is worth mentioning since one of the main paradigms of geometric integration is to develop numerical schemes with good qualitative behavior in long-time simulations. One potential drawback of the proposed integrator is the execution time which is significantly larger than the ode15i. However, we point out that this is mainly because the ode15i is an optimized package designed for Matlab, while the implementation of the proposed integrator makes use of a particular subroutine implementing a modified multivariate Newton-Raphson method, which has not been optimized for the Matlab platform.

It is also worth noticing that the anomaly displayed by the Lagrange-d'Alembert integrator (regarding Figure \ref{fig:foucault_sim_2_oscillation_plane}) disappears when we refine the time-step by a factor of 20, which certainly implies a computational cost that might be undesirable. This raises the question of whether this anomaly is geometric in nature or is it a purely numerical issue.

\section*{Acknowledgments}
C.E.S. acknowledges the support given by PRONII. I.O. and E.M. acknowledge partial support given by  FEEI-CONACYT-PROCIENCIA.

\clearpage
\bibliographystyle{siamplain}
\bibliography{references}

\appendix \section{Discrete equations of the contact integrator with constraints for the falling disk}

The following equations enumerated from \eqref{eq:first_equation_of_the_integrator} to (\ref{eq:last_equation_of_the_integrator}) correspond to the resulting equations of the integrator based on the discretization of the Herglotz variational principle with constraints that we propose in this work, applied to the problem of the falling disk described in Section \ref{sec:falling_disk}.
\small

\begin{equation}\label{eq:first_equation_of_the_integrator}
    \begin{split}
      &R \cos(\phi_j) \left[\frac{F^X_j}{2} - \frac{1}{\frac{\alpha h}{2} + 1} \left(\frac{F^X_j}{2} - \frac{m (X_{j-1} - X_j)}{h^2}\right) \left(\frac{\alpha h}{2} - 1\right) + \frac{m (X_j - X_{j+1})}{h^2} \right] \\
  &+ R \sin(\phi_j) \left[ \frac{F^Y_j}{2} - \frac{1}{\frac{\alpha h}{2} + 1} \left(\frac{F^Y_j}{2} - \frac{m (Y_{j-1} - Y_j)}{h^2}\right) \left(\frac{\alpha h}{2} - 1\right) + 
  \frac{m (Y_j - Y_{j+1})}{h^2} \right]\\
  &-\frac{1}{\frac{\alpha h}{2} + 1}
  \left(\frac{\alpha h}{2} - 1\right) 
  \left[\frac{F^{\psi}_j}{2} - \frac{1}{h} I_A \left( \frac{\psi_{j-1} - \psi_j}{h} - \sin\left(\frac{\theta_{j-1} + \theta_j}{2}\right) \frac{\phi_{j-1} - \phi_j}{h}\right)\right] \\
  &+\frac{1}{h} I_A \left( \frac{\psi_j - \psi_{j+1}}{h} - 
  \sin\left(\frac{\theta_j + \theta_{j+1}}{2}\right) 
  \frac{\phi_j - \phi_{j+1}}{h} \right) +\frac{F^{\psi}_j}{2}= 0.
    \end{split}
\end{equation}


\begin{equation}
    \begin{split}
      &R \cos\left(\frac{\phi_j + \phi_{j+1}}{2}\right) \frac{\psi_j - \psi_{j+1}}{h} - 
  \frac{X_j - X_{j+1}}{h}\\
  &- R \cos\left(\frac{\phi_j + \phi_{j+1}}{2}\right) \sin\left(\frac{\theta_j + \theta_{j+1}}{2}\right) \frac{\phi_j - \phi_{j+1}}{h}\\
  &-R \cos\left(\frac{\theta_j + \theta_{j+1}}{2}\right) \sin\left(\frac{\phi_j + \phi_{j+1}}{2}\right) 
  \frac{\theta_j - \theta_{j+1}}{h} = 0.
    \end{split}
\end{equation}


\begin{equation}
    \begin{split}
      &R \sin\left(\frac{\phi_j + \phi_{j+1}}{2}\right) \frac{\psi_j - \psi_{j+1}}{h} - \frac{Y_j - Y_{j+1}}{h}\\ 
      &- R \sin\left(\frac{\phi_j + \phi_{j+1}}{2}\right) \sin\left(\frac{\theta_j + \theta_{j+1}}{2}\right) 
  \frac{\phi_j - \phi_{j+1}}{h}\\
  &+R \cos\left(\frac{\phi_j + \phi_{j+1}}{2}\right) \cos\left(\frac{\theta_j + \theta_{j+1}}{2}\right) \frac{\theta_j - \theta_{j+1}}{h} = 0.
    \end{split}
\end{equation}

\begin{equation}
    \begin{split}
        &\frac{F^{\theta}_j}{m} + \frac{I_T}{m} \left(\frac{2(\theta_j - \theta_{j+1})}{h^2} - \frac{\sin(\theta_j+\theta_{j+1})}{2} \frac{(\phi_j - \phi_{j+1})^2}{h^2}\right)+R g\sin(\theta_j)\\
  &+R^2\left( \sin^2\left(\frac{\theta_j + \theta_{j+1}}{2}\right) \frac{2(\theta_j - \theta_{j+1})}{h^2} + \frac{\sin(\theta_j+\theta_{j+1})}{2} \frac{(\theta_j - \theta_{j+1})^2}{h^2}\right)\\
  &+\frac{\alpha h-2}{\alpha h+2} \Bigg\{ \frac{I_T}{m} \left[\frac{2(\theta_{j-1} - \theta_j)}{h^2} + \frac{\sin(\theta_{j-1}+\theta_j)}{2}\frac{(\phi_{j-1} - \phi_j)^2}{h^2}\right] - \frac{F^{\theta}_j}{m}- R g \sin(\theta_j)\\
  &+R^2 \left[\sin^2\left(\frac{\theta_{j-1} + \theta_j}{2}\right) \frac{2(\theta_{j-1} -  \theta_j)}{h^2} - \frac{\sin(\theta_{j-1}+\theta_j)}{2} \frac{(\theta_{j-1} - \theta_j)^2}{h^2}\right]\\
  &+\frac{I_A (\phi_{j-1} - \phi_j)}{m} \cos\left(\frac{\theta_{j-1} + \theta_j}{2}\right) \left[\frac{\psi_{j-1} - \psi_j}{h^2} - \sin\left(\frac{\theta_{j-1} + \theta_j}{2}\right) \frac{\phi_{j-1} - \phi_j}{h^2}\right] \Bigg\}\\
  &+2R \cos(\phi_j) \cos(\theta_j) \left[\frac{F^Y_j}{2m} - \frac{\alpha h-2}{\alpha h+2}\left(\frac{F^Y_j}{2m} - \frac{Y_{j-1} - Y_j}{h^2}\right) + \frac{Y_j - Y_{j+1}}{h^2}\right]\\
  &-2R \cos(\theta_j) \sin(\phi_j) \left[\frac{F^X_j}{2m} - \frac{\alpha h-2}{\alpha h+2}\left(\frac{F^X_j}{2m} - \frac{X_{j-1} - X_j}{h^2}\right) + \frac{X_j - X_{j+1}}{h^2} \right]\\
  & - \frac{I_A(\phi_j - \phi_{j+1})}{m}\cos\left(\frac{\theta_j + \theta_{j+1}}{2}\right)  \left[\frac{\psi_j - \psi_{j+1}}{h^2} - \sin\left(\frac{\theta_j + \theta_{j+1}}{2}\right) \frac{\phi_j - \phi_{j+1}}{h^2}\right]  = 0.
    \end{split}
\end{equation}


\begin{equation}\label{eq:last_equation_of_the_integrator}
    \begin{split}
        &\frac{F^{\phi}_j}{2}+ I_T \cos^2\left(\frac{\theta_j + \theta_{j+1}}{2}\right) \frac{\phi_j - \phi_{j+1}}{h^2}\\
        &- \frac{\alpha h - 2}{\alpha h+2} \Bigg\{ \frac{F^{\phi}_j}{2} - I_T \cos^2\left(\frac{\theta_{j-1} + \theta_j}{2}\right) \frac{\phi_{j-1} - \phi_j}{h^2}\\
  &+I_A \sin\left(\frac{\theta_{j-1} + \theta_j}{2}\right) \left[\frac{\psi_{j-1} - \psi_j}{h^2} - \sin\left(\frac{\theta_{j-1} + \theta_j}{2}\right) \frac{\phi_{j-1} - \phi_j}{h^2}\right] \Bigg\}\\
 &- R \cos(\phi_j) \sin(\theta_j) \left[\frac{F^X_j}{2} - \frac{\alpha h - 2}{\alpha h+2} \left(\frac{F^X_j}{2} - \frac{m (X_{j-1} - X_j)}{h^2}\right) + \frac{m (X_j - X_{j+1})}{h^2}\right]\\
 &- R \sin(\phi_j) \sin(\theta_j) \left[\frac{F^Y_j}{2} - \frac{\alpha h - 2}{\alpha h+2} \left(\frac{F^Y_j}{2} - \frac{m (Y_{j-1} - Y_j)}{h^2}\right) + \frac{m (Y_j - Y_{j+1})}{h^2}\right]\\
 & -
  I_A \sin\left(\frac{\theta_j + \theta_{j+1}}{2}\right) \left[\frac{\psi_j - \psi_{j+1}}{h^2} -
  \sin\left(\frac{\theta_j + \theta_{j+1}}{2}\right) \frac{\phi_j - \phi_{j+1}}{h^2}\right]= 0.
    \end{split}
\end{equation}

\noindent {\it E-mail address:} \email{emaciel@pol.una.py}\\
{\it E-mail address:} \email{iortiz@pol.una.py}\\
{\it E-mail address:} \email{cschaer@pol.una.py}

\end{document}

%% file: ex_shared.tex
\usepackage{lipsum}
\usepackage{amsfonts}
\usepackage{graphicx}
\usepackage{epstopdf}
\usepackage{algorithmic}
\ifpdf
  \DeclareGraphicsExtensions{.eps,.pdf,.png,.jpg}
\else
  \DeclareGraphicsExtensions{.eps}
\fi

\newsiamremark{remark}{Remark}
\newsiamremark{hypothesis}{Hypothesis}
\crefname{hypothesis}{Hypothesis}{Hypotheses}
\newsiamthm{claim}{Claim}
\newsiamthm{example}{Example}

\headers{Nonholonomic Herglotz Integrator}{E. Maciel, I. Ortiz, and C. E. Schaerer}

\title{A Herglotz-based integrator for nonholonomic mechanical systems}

\author{Elias Maciel\and Inocencio Ortiz\and Christian E. Schaerer}
